\documentclass[10pt]{article}
\usepackage[OT2,T1]{fontenc}
\usepackage{amssymb,amsfonts,amsmath,amsthm,amscd}
\usepackage[margin=1in]{geometry}
\usepackage{authblk}
\usepackage[inline]{enumitem}
\usepackage{graphicx}
\usepackage{float}
\usepackage{cite}

\def\TT{\mathsf{T}}
\def\tt{\mathfrak{t}}
\def\HH{\mathsf{H}}
\def\hh{\mathfrak{h}}
\def\MM{\mathsf{M}}
\def\DD{\mathcal{D}}
\def\RR{\mathbb{R}}
\def\CC{\mathbb{C}}
\def\SL{\rm{SL}_2(\CC)}
\DeclareMathOperator{\diag}{diag}
\newtheorem{theorem}{Theorem}
\newtheorem{lemma}{Lemma}
\newtheorem{proposition}{Proposition}
\newtheorem{corollary}[theorem]{Corollary}
\theoremstyle{definition}

\newtheorem{remark}[theorem]{Remark}
\newtheorem*{thank}{{\indent Acknowledgement}}
\usepackage{xcolor}
\usepackage{verbatim}

\title{Integrable and non-integrable LR geodesic flows on 5-dimensional 2-step nilpotent Lie groups}

\author{Milan Pavlovi\'{c}, Tijana \v{S}ukilovi\'{c}}
\affil{University of Belgrade, Faculty of Mathematics, \\Studentski Trg 16 11000 Belgrade, Serbia\\ milan.pavlovic@matf.bg.ac.rs, tijana.sukilovic@matf.bg.ac.rs}

\allowdisplaybreaks

\begin{document}

\maketitle

\begin{abstract}
We investigate sub-Riemannian geodesic flows on 5-dimensional nilpotent Lie groups endowed with left-invariant Riemannian metrics. Up to isomorphism, there are exactly two such groups: the Heisenberg group $H_5$ and the 5-dimensional nilpotent Lie group with a 2-dimensional commutator subgroup. For each group, we study geodesic flows associated with both left-invariant (LL) and right-invariant (LR) distributions, deriving explicit descriptions of the corresponding geodesics and comparing their integrability properties. While the LL systems are known to be Liouville integrable, and even superintegrable, we show that LR systems are Liouville integrable only in two cases: on the Heisenberg group and, for the second group, on a distinguished invariant level set.


\end{abstract}

\section{Introduction}

Sub-Riemannian geometry extends Riemannian geometry by imposing constraints on the admissible directions of motion, thereby highlighting its inherently anisotropic nature. Within this framework, nilpotent groups constitute a particularly important class of spaces. They arise as tangent cones of sub-Riemannian manifolds, playing a role analogous to that of tangent spaces in Riemannian geometry. 

The most extensively studied sub-Riemannian structures on Lie groups are those defined by a left-invariant Riemannian metric together with a left-invariant distribution. Because left translations are sub-Riemannian isometries that act transitively on the group and preserve geodesics, the problem of joining arbitrary points $P$ and $Q$  reduces to joining the identity element with $P^{-1}Q$. We refer to these structures as \textit{LL systems} (for more details see~\cite{ABB19,BN13,M02}). 

In contrast, \textit{LR systems} are characterized by a left-invariant Riemannian metric and a right-invariant distribution. The fact that the kinetic energy and the velocity constraints are governed by different symmetry groups gives rise to more sophisticated and intricate geometry. Such systems appear naturally in nonholonomic mechanics~\cite{VV88,FJ04} and have also been investigated from the perspective of sub-Riemannian  geometry~\cite{T97,AB17,PS}.

In this paper, we investigate the integrability of sub-Riemannian geodesic flows on 5-dimensional nilpotent Lie groups equipped with a left-invariant Riemannian metric. There are exactly two 5-dimensional nilpotent Lie groups of this type, namely the Heisenberg group $H_5$ and the group $\TT$ with a 2-dimensional commutator subgroup. Our focus is on both classes of systems introduced above, namely LL systems and LR systems. Although the integrability of LL systems has been studied previously, we include a concise and unified treatment to provide a convenient framework for comparison with LR systems, whose integrability properties are investigated here for the first time. The main objective is to determine and compare their integrability properties, thereby highlighting the influence of the underlying symmetry of the distribution on the dynamics of the geodesic flow. This is motivated by the paper~\cite{M97} in which the authors ask whether the geodesic flow on a Carnot group is always integrable and establishes an affirmative answer for LL systems on 2-step nilpotent groups (the proof of this in 3-dimensional case can be found in~\cite{ABB19}), but provide an example of a 3-step Carnot group whose sub-Riemannian geodesic flow is not algebraically completely integrable. In the present context, we consider an analogous question for LR systems on 2-step nilpotent groups, showing that integrability holds only under certain special conditions and that, in general, there exist examples of non-integrable systems.

Before proceeding to the individual cases, we formulate the main results of the paper. The proofs are deferred to later sections, where the sub-Riemannian structures are analyzed separately and the necessary computations are carried out in detail.

The integrability of sub-Riemannian LL systems considered in this paper was established in~\cite{MPAM99,BBG12} and is summarized in the following theorem.
\begin{theorem}\label{thm:LL_main1}
    For a 5-dimensional 2-step nilpotent Lie group $G$ endowed with a left-invariant sub-Riemannian structure $(\DD, g)$, the associated geodesic flow defines a Hamiltonian system that is Liouville integrable and, in fact, superintegrable.
\end{theorem}

While it is a classical fact that the normal geodesics of the Heisenberg group $H_5$ are horizontal lifts of segments of torus windings in $\RR^4$  (sometimes called sub-Riemannian helices), it is noteworthy that, in the case of the group $\TT$, the geodesics themselves are generalized helices in the usual sense. Finally, we recall that, since the distributions are 2-step bracket generating, every abnormal length minimizer is also normal, and hence there are no strict abnormal geodesics (see~\cite{BBG12,AS04}).

\begin{theorem}\label{thm:LL_main2}
    The normal geodesics associated with the left-invariant sub-Riemannian structure $(\DD, g)$ on a 5-dimensional 2-step nilpotent Lie group $G$ are horizontal lifts of helix arcs; straight lines are recovered as a limiting case. When $G=\TT$, the resulting geodesics are generalized helices in the classical sense.
\end{theorem}

Whereas LL systems are superintegrable, LR systems are Liouville integrable only in two cases: when the underlying group is $H_5$, or when the dynamics is reduced to a special level set.

\begin{theorem}\label{thm:LR_main1}
    Among all sub-Riemannian LR Hamiltonian systems arising from 5-dimensional 2-step nilpotent Lie groups, the system associated with the Heisenberg group $H_5$ is the only integrable one. All remaining systems are non-integrable.
\end{theorem}
\begin{proof}
The integrability of the Hamiltonian system on the Heisenberg group $H_5$, associated with a sub-Riemannian structure defined by a left-invariant Riemannian metric and a right-invariant distribution, is established in  Corollary~\ref{cor:HLR}. 
The LR systems defined on the group $\TT$ are analyzed in Theorems~\ref{thm:LR3} and~\ref{th:LR4}, which prove their non-integrability.
\end{proof}

\begin{theorem}\label{thm:LR_main2}
    The reduced Hamiltonian system obtained by symplectic reduction of the LR system is equivalent to a Hamiltonian system describing the motion of a charged particle in the constant magnetic field. The dynamics on certain special level sets are integrable.
\end{theorem}
\begin{proof}
Propositions~\ref{prop:magnetic1} and~\ref{prop:magnetic2}, and Theorem~\ref{thm:rightR} highlight the magnetic nature of the systems under consideration. For the Heisenberg group $H_5$, the reduction of the dynamics to any level set $\lambda_5=C$ yields an integrable system (see Theorem~\ref{thm:rightR}). For the group $\TT$, Proposition~\ref{prop:heis} provides an example in which the reduced dynamics is equivalent to the LR Hamiltonian system on the Heisenberg group $H_3$.
\end{proof}

The remainder of the paper is organized as follows. In Section~\ref{sec:pre} we briefly recall the necessary background on sub-Riemannian geometry and integrability. Sections~\ref{sec:T} and~\ref{sec:H}  are devoted to the two nilpotent Lie groups under consideration: first to the group $\TT$ and the second to $H_5$. 
For each group, we study both the LL and LR sub-Riemannian structures and analyse the integrability of the corresponding geodesic flows.

\section{Preliminaries}\label{sec:pre}

Recall that a rank-$k$ distribution $\DD \subset T\mathsf{M}$ on a connected $n$-dimensional manifold $\mathsf{M}$ is called completely nonholonomic (or bracket-generating) if the space $\Gamma(\DD)$ of its sections generates the Lie algebra $\mathcal{X}(\mathsf{M})$ of vector fields. Locally, the distribution $\DD$ can be described as the kernel of a system $(\omega^1, \cdots, \omega^{n-k})$ of $1$-forms. A sub-Riemannian manifold $(\mathsf{M}, \DD, g)$ consists of a manifold $\mathsf{M}$, a distribution $\DD$, and a Riemannian metric $g$ defined on $\DD$.

With reference to~\cite{M02}, we present the following theorems without proof.

\begin{theorem}[on normal geodesics]\label{thm:normal_geodesics} Let $\xi (t)= (\gamma (t), \lambda (t))$ be a solution of the Hamilton's differential equations on $T^*Q$ for the sub-Riemannian Hamiltonian $H$ and let $\gamma (t)$ be its projection on $Q$. Then every sufficiently short arc of $\gamma$ is a minimizing sub-Riemannian geodesic. Moreover, $\gamma$ is the unique minimizing geodesic connecting its endpoints.
\end{theorem}

The projected curves $\gamma(t)$ in Theorem~\ref{thm:normal_geodesics} are called \textit{normal} sub-Riemannian geodesics. Minimizing geodesics that do not arise as solutions of the Hamiltonian system are called \textit{singular} or \textit{abnormal} geodesics.

\begin{theorem}[Chow–Rashevskii]\label{thm:Chow–Rashevskii}
 If $\mathcal{H}$ is a bracket-generating distribution on a connected manifold $\MM$ then any two points of $\MM$ can be joined by a horizontal path.
\end{theorem}

Hamiltonian flow is given by the Hamiltonian function. 
\begin{align*}
H(\mathbf{x},\mathbf{\lambda} )=\frac{1}{2}\sum_{i,j=1}^{5}g^{ij} \lambda _i \lambda _j,
\end{align*}
and normal geodesic corresponds to the solution of the Hamilton system:
\begin{align}\label{Hamilton_system}
  \dot{x}_i &= \frac{\partial H}{\partial \lambda_i}, \qquad
  \dot{\lambda}_i = -\frac{\partial H}{\partial x_i},
\end{align}
where $(\mathbf{x},\mathbf{\lambda})=(x_i,\lambda_i)_{1\leq i \leq n}$ is the point in the cotangent bundle of $\mathsf{M}$.

If the horizontal vector fields spanning the distribution $\DD$ at a point $p$, expressed in local coordinates, are given by
\begin{align*}
X_j=\sum_{i=1}^{n}x_i^j(p)\frac{\partial}{\partial x_i}, \quad j=1,\ldots k,
\end{align*}
then, as shown in~\cite{CC09}, the coefficients of the quadratic form defining the Hamiltonian function $H$ are given~by
\begin{align}\label{coef:HamiltonsEq}
g^{ij}=\sum_{l=1}^{k}x_i^l x_j^l, \quad 1\leq i, j\leq n.
\end{align}

Let $f$ and $g$ be first integrals of the Hamiltonian system~\eqref{Hamilton_system} with Hamiltonian $H$, that is,
$\{f,H\}=\{g,H\}=0$, where $\{\cdot,\cdot\}$ denotes the standard Poisson bracket. By the Jacobi identity, the Poisson bracket $\{f,g\}$ is itself a first integral. Consequently, the collection of all first integrals forms a Lie algebra $\mathcal{F}$ with respect to the Poisson bracket.

A Hamiltonian system on a $2n$-dimensional phase space~\eqref{Hamilton_system} is called \emph{completely integrable in the non-commutative sense} if it admits a complete algebra $\mathcal{F}$ of first integrals. According to a conjecture of  Mishchenko and Fomenko~\cite{MF78}, non-commutative integrability should imply Liouville integrability by means of first integrals that belong to the same functional class as those in $\mathcal{F}$ and commute pairwise with respect to the Poisson bracket. The classical Liouville integrability is recovered when $\mathcal{F}$ itself is a complete commutative algebra. If, in addition, the system possesses more than $n$ functionally independent first integrals, it is said to be \emph{superintegrable}. In the  $C^\infty$-smooth category, non-commutative integrability implies Liouville integrability, as proven in~\cite{BJ03}.

Non-integrability can be established using differential Galois theory. According to~\cite{MRR01}, Liouville integrability imposes strong restrictions on the differential Galois group of the normal variational equation along any non-equilibrium solution: its identity component must be abelian. Consequently, the presence of a non-abelian differential Galois group provides an obstruction to Liouville integrability.
\begin{theorem}[Morales–Ramis]\label{thm:MR}
    If a Hamiltonian system is Liouville integrable, with meromorphic first integrals, then the identity component of the differential Galois group of variational equations along any non-constant solution is abelian.
\end{theorem}
This criterion will serve as the main tool in our proofs of non-integrability throughout the paper.

\section{Group $\TT$}\label{sec:T}

There exists a unique 5-dimensional, 2-step nilpotent, simply connected real Lie group with a 2-dimensional commutator subgroup. Following~\cite{BN23}, we denote this group by $\TT$:
\begin{align*}
\TT = \left\{
\begin{bmatrix}
1 & x_1 & x_4 & x_5 \\
0 & 1 & x_2 & x_3 \\
0 & 0 & 1 & 0 \\
0 & 0 & 0 & 1 
\end{bmatrix}: x_1, x_2, x_3, x_4, x_5 \in \mathbb{R} \right\}.
\end{align*}
The group $\TT$ is defined with respect to matrix multiplication. Its corresponding Lie algebra $\tt$ is identified with the tangent space $T_e\TT$ at the identity element $e$, and there exists a basis ${e_i}_{1 \leq i \leq 5}$ in which the commutation relations are given by:
\begin{align*}
[e_1,e_2]=e_4, \quad [e_1,e_3]=e_5.
\end{align*}

In~\cite{BN23} authors have done full classification of left-invariant sub-Riemannian structures and it is proved that, in addition to the structure obtained from the whole Lie algebra $\tt=Lie(\TT)$ which is actually Riemannian, there are only two non-isometric  left-invariant sub-Riemannian structures.

\begin{lemma}[\!\!\cite{BN23}]
    Any left-invariant sub-Riemannian structure $(\DD, g)$ on $\TT$ is isometric to exactly one of the structures:
        \begin{align*}
            (\TT,\mathcal{L}_3,g): &\quad \begin{cases}
                \mathcal{L}_3(e)= \RR\langle e_1, e_2, e_3\rangle\\
                g=\diag(1,1,1)
            \end{cases}\\
            (\TT,\mathcal{L}_4,h): &\quad \begin{cases}
                \mathcal{L}_4(e)= \RR\langle e_1, e_2, e_3, e_4\rangle\\
                h=\mu\cdot \diag(1,1,1,1),\ \mu>0.
            \end{cases}
    \end{align*}
\end{lemma}

Neither of these distributions is contact. Let 
$\alpha=dx_4-x_1dx_2$ and $\beta = dx_5-x_1dx_3$ be $1-$forms, then $\mathcal{L}_3=\ker \alpha \cap \ker\beta$ and $\mathcal{L}_4=\ker\beta$. The distribution $\mathcal{L}_3$ cannot be contact, as its rank is 3. Moreover, $\mathcal{L}_4$ is not contact, since $\beta\wedge(d\beta)^2=0$.

In the sequel, we study the LL and LR systems corresponding to each of these two sub-Riemannian structures, and, for each system, derive and visualize the corresponding normal sub-Riemannian geodesics.

\subsection{LL systems}

The group $\TT$ acts on itself by left translations $L_g : \TT \to \TT$, defined by $L_g(h) = gh$ for $g, h \in \TT$. For $g = (x_1, x_2, x_3, x_4, x_5) \in \TT$, the differential of this action is given by:
\begin{align*}
dL_g(e_1)=e_1, \quad dL_g(e_2)=e_2+x_1 e_4, \quad dL_g(e_3)=e_3 + x_1 e_5, \quad dL_g(e_4)=e_4, \quad dL_g(e_5)=e_5. 
\end{align*}

The left-invariant distributions generated by the subspaces $\mathcal{L}_3 = \mathbb{R}\langle e_1, e_2, e_3\rangle$ and $\mathcal{L}_4 = \mathbb{R}\langle e_1, e_2, e_3, e_4\rangle$ are, respectively, the 3-dimensional distribution $\DD^3 = dL_g \mathcal{L}_3$ and the 4-dimensional distribution $\DD^4 = dL_g \mathcal{L}_4$. It is immediate that both distributions are completely nonholonomic. 

While~\cite{BBG12} establishes the integrability of sub-Riemannian LL-geodesics in general 2-step Carnot groups of corank $2$ and provides a framework for their computation, here we compute them explicitly and derive several of their geometric properties.

\subsubsection{Distribution $\DD^3$}

A sub-Riemannian metric $g^{ij}$ on the distribution $\DD^3$ is given by  $g^{ij}(e) = \diag(1,1,1,0,0)$ or in local coordinates:
\begin{align*}
g^{ij}(p) = \begin{bmatrix}
1 & 0 & 0 & 0 & 0 \\
0 & 1 & 0 & x_1 & 0 \\
0 & 0 & 1 & 0 & x_1 \\
0 & x_1 & 0 & x_1^2 & 0 \\
0 & 0 & x_1 & 0 & x_1^2
\end{bmatrix}.
\end{align*}

This sub-Riemannian metric is compatible with the standard left-invariant Riemannian metric $\langle \cdot, \cdot \rangle$ on $\TT$, whose coordinate representation is given by a symmetric matrix:
\begin{align}\label{eq:S}
S= \begin{bmatrix}
1 & 0       & 0       & 0    & 0 \\
0 & 1+x_1^2 & 0       & -x_1 & 0 \\
0 & 0       & 1+x_1^2 & 0    & -x_1 \\
0 & -x_1    & 0       & 1    & 0 \\
0 & 0       & -x_1    & 0    & 1
\end{bmatrix}. 
\end{align}

In this case, a direct computation yields the Hamiltonian
\begin{align*}
H(\mathbf{x},\mathbf{\lambda} )=\frac{1}{2}\sum_{i,j=1}^{5}g^{ij} \lambda _i \lambda _j=\frac{1}{2}(\lambda_1^2+(\lambda_2+x_1\lambda_4)^2+(\lambda_3+x_1\lambda_5)^2),
\end{align*}
and the associated Hamiltonian system:
\begin{gather} \label{system:Ham3D}
  \begin{aligned}
  \dot{x}_1 &= \lambda_1, &\dot{\lambda} _1 &= -\lambda_2\lambda_4-\lambda_3\lambda_5-x_1(\lambda_4^2+\lambda_5^2),\\
  \dot{x}_2 &= \lambda_2+x_1\lambda_4, \qquad\qquad & \dot{\lambda} _2 &= 0,\\
  \dot{x}_3 &= \lambda_3+x_1\lambda_5,  & \dot{\lambda} _3 &= 0,\\
  \dot{x}_4 &= x_1(\lambda_2+x_1 \lambda_4), &  \dot{\lambda} _4 &= 0,\\
  \dot{x}_5 &= x_1(\lambda_3+x_1\lambda_5),   & \dot{\lambda} _5 &= 0.
 \end{aligned}
\end{gather}
It follows immediately that this Hamiltonian system  possesses six functionally independent first integrals:
\begin{align}\label{eq:integrali}
I_0 = H,\quad I_1 = \lambda_1+x_2\lambda_4+x_3\lambda_5,\quad I_2 =\lambda_2, \quad I_3 =\lambda_3, \quad
I_4 =\lambda_4, \quad I_5 =\lambda_5.
\end{align}
The integrals $I_0, I_2,\ldots, I_5$ are pairwise in involution, making the system completely integrable in Liouville sense. Additionally, the existence of the sixth integral $I_1$ such that $I_4 =\{I_1, I_2\}$ and $ I_5 =\{I_1, I_3\}$ makes the system superintegrable.

Having determined the first integrals of the system, we introduce the substitution 
\begin{align*}
    \lambda_4 = \omega \cos \phi,\quad \lambda_5 = \omega \sin \phi, \quad \omega \ge 0,\ \phi \in [0, 2\pi). 
\end{align*}
This allows us to characterize the geodesic lines passing through the origin. Since the proof consists of direct calculations, we omit the details and state the result.

\begin{theorem}\label{geodesics:D3}
For any point $P \in \TT$, there exists a normal sub-Riemannian geodesic $\gamma$ for the sub-Riemannian structure $(\TT, \DD^3, g)$ joining $P$ to the origin $O = (0,0,0,0,0)$, satisfying $\gamma(0)=\mathbf{0}$, and whose expression in local coordinates is given by:
\begin{enumerate}[wide, topsep=0pt, itemsep=1pt, labelindent=10pt, label={Case}]
\item $\omega = 0$:
$\gamma(t)=(\lambda_1 t, \lambda_2 t, \lambda_3 t, \frac{1}{2}t^2\lambda_1\lambda_2,\frac{1}{2}t^2\lambda_1\lambda_3)$,
where $\lambda_1,\lambda_2,\lambda_3\in\mathbb{R}.$
\item $\omega \neq 0$: $\gamma(t)=(x_1(t),x_2(t),x_3(t),x_4(t),x_5(t))$ where $c,\lambda_2,\lambda_3\in\mathbb{R}$ and the coordinate functions $x_i(t), i\in\{1,\ldots,5\}$  are as follows:
 \begin{align*} 
x_1(t)&= c \sin (t \omega ){-}\frac{2 \sin ^2\left(\frac{t \omega }{2}\right) \left(\lambda _3 \sin (\phi
   ){+}\lambda _2 \cos (\phi )\right)}{\omega }\\
x_2(t)&=\cos (\phi ) \left(2c\sin^2\left(\frac{t\omega}{2}\right){+}\frac{\sin (t \omega ) \left(\lambda _2 \cos (\phi
   ){+}\lambda _3 \sin(\phi )\right)}{\omega }\right){+}t \sin (\phi ) \left(\lambda _2
   \sin (\phi ){-}\lambda _3 \cos (\phi )\right) \\
x_3(t)&=\sin (\phi )\left(2c\sin^2\left(\frac{t\omega}{2}\right){+}\frac{ \sin (t \omega )\left(\lambda _2 \cos (\phi
   ){+}\lambda _3 \sin(\phi ) \right)}{\omega }\right){+}t \cos (\phi ) \left(\lambda _3 \cos (\phi
   ){-}\lambda _2 \sin (\phi )\right)\\
x_4(t)&= {-}\frac{c \cos (\phi ) \cos (2 t \omega ) \left(\lambda _3 \sin (\phi ){+}\lambda _2 \cos (\phi)\right)}{2 \omega }{+}\frac{c \cos (t \omega ) \left(\lambda _3 \sin (2 \phi ){+}\lambda _2 \cos
   (2 \phi )\right)}{\omega }\\
   &{+}\frac{\left(3 \lambda _3 \sin (2 \phi ){+}\lambda _2 (3 \cos (2 \phi
   ){-}1)\right) \left({-}c{+}t(\lambda _3 \sin (\phi ){+}\lambda _2  \cos (\phi )\right)}{4 \omega
   }{+}\frac{ c^2}{4} \cos (\phi ) (2 t \omega {-}\sin (2 t \omega ))\\
   &{-}\frac{\sin (t \omega )
   \left(2 \lambda _2 \lambda _3 \sin (3 \phi ){+}\lambda _2^2 (\cos (\phi ){+}\cos (3 \phi )){+}4 \lambda
   _3^2 \sin ^2(\phi ) \cos (\phi )\right)}{2 \omega ^2}\\
   &{+}\frac{\cos (\phi ) \sin (2 t \omega )
   \left(\lambda _3 \sin (\phi ){+}\lambda _2 \cos (\phi )\right){}^2}{4 \omega ^2}\\
x_5(t)&=\frac{c \cos (t \omega ) \left(\lambda _2 \sin (2 \phi ){-}\lambda _3 \cos (2 \phi )\right)}{\omega
   }{-}\frac{c \sin (\phi ) \cos (2 t \omega ) \left(\lambda _3 \sin (\phi ){+}\lambda _2 \cos (\phi
   )\right)}{2 \omega }\\
   &{+}\frac{\left(3 \lambda _2 \sin (2 \phi ){-}\lambda _3 (3 \cos (2 \phi
   ){+}1)\right) \left({-}c{+}t(\lambda _3 \sin (\phi ){+}\lambda _2 \cos (\phi ))\right)}{4 \omega
   }{+}\frac{c^2}{4}  \sin (\phi ) (2 t \omega {-}\sin (2 t \omega ))\\
   &{+}\frac{\sin (t \omega )
   \left(\lambda _3^2 (\sin (3 \phi ){-}\sin (\phi )){+}2 \lambda _3 \lambda _2 \cos (3 \phi ){-}4 \lambda
   _2^2 \sin (\phi ) \cos ^2(\phi )\right)}{2 \omega ^2}\\
   &{+}\frac{\sin (\phi ) \sin (2 t \omega )
   \left(\lambda _3 \sin (\phi ){+}\lambda _2 \cos (\phi )\right){}^2}{4 \omega ^2}.
\end{align*} 
\end{enumerate}
\end{theorem}

We investigate more closely the projection of $\gamma$ on a 3-dimensional space, $\gamma_{123}(t)= \pi_{123}(\gamma(t)) = (x_1(t), x_2(t), x_3(t))$.
\begin{theorem}\label{thm: helices_D3}
Let $\gamma:\ \mathbb{R}\rightarrow \TT$ be normal sub-Riemannian geodesic  and $\gamma_{123}:\ \mathbb{R}\rightarrow {\mathbb{R}^3}$ its projection on left-invariant distribution generated with $\left\{X_1,X_2,X_3 \right\}$, then $\gamma$ is a generalized helix in dimension $5$ and $\gamma_{123}$ is a curve of constant curvature and constant torsion, i.e. circular helix.
\end{theorem}

\begin{proof}
To prove that $\gamma$ is generalized helix it is enough to take $v=(0,-\sin{\phi},\cos{\phi},0,0)$ and calculate 
\begin{align*}
    \langle\gamma'(t),v\rangle=-\lambda_2\sin{\phi}+\lambda_3\cos{\phi}\equiv const=\kappa.
\end{align*} 

For $\gamma_{123}$, provided it is arc length parametrized, i.e. $\lambda_2^2+\lambda_3^2+c^2\omega^2=1$, we calculate curvature and torsion as
\begin{flalign*}
&&k(t)&=\omega\sqrt{(1+\lambda_2\sin\phi-\lambda_3\cos\phi)(1-\lambda_2\sin\phi+\lambda_3\cos\phi)}=\omega \sqrt{1-\kappa^2}&\\
\text{and}&&
\tau(t)&=\omega\kappa, &
\end{flalign*}
which are constants.
\end{proof}
\begin{figure}[H]
\begin{center}
\includegraphics[width=0.45\linewidth]{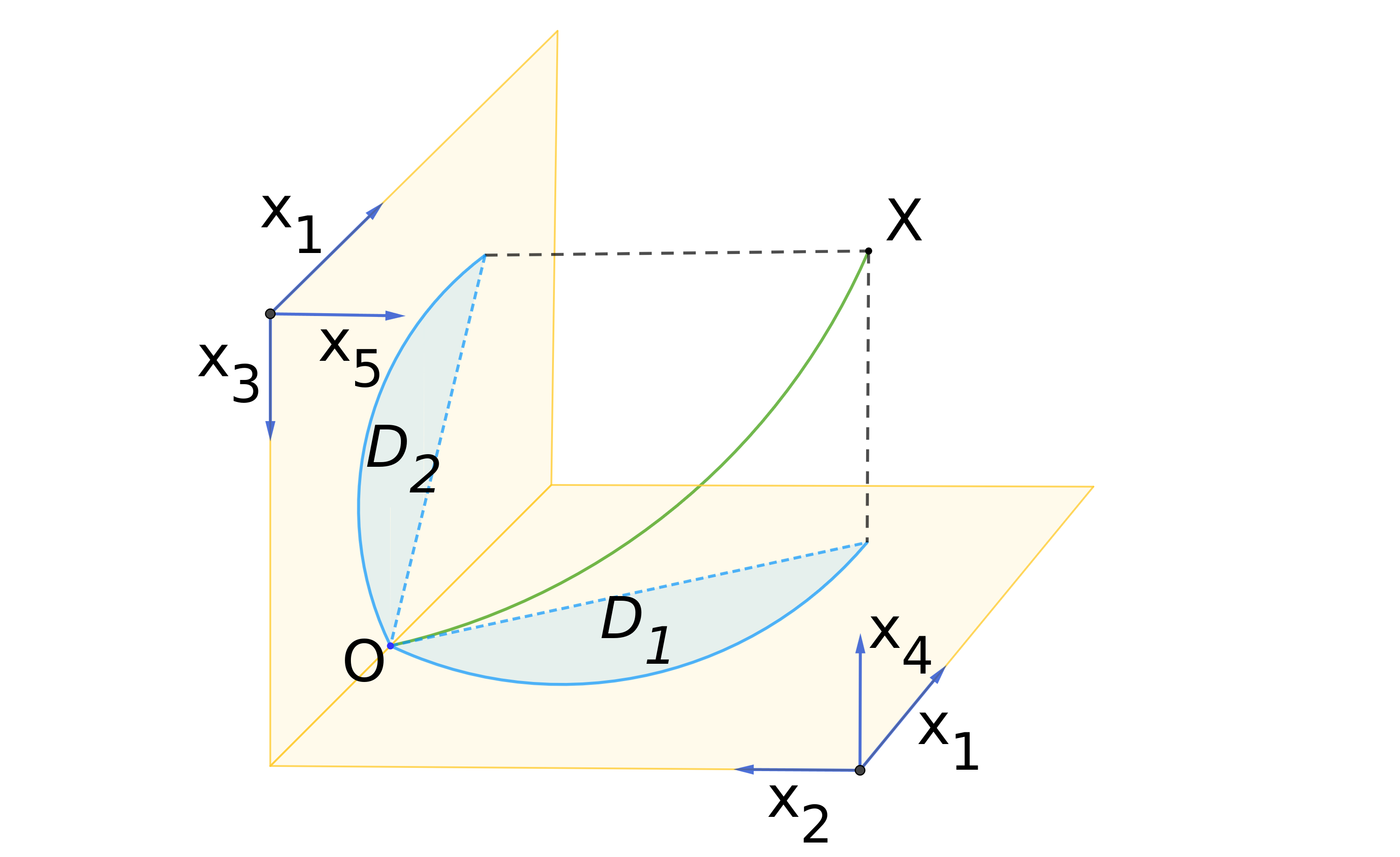}
\caption{Projections of a sub-Riemannian geodesic onto the $x_1x_2$- and $x_1x_3$-planes.}\label{fig:ds}
\end{center}
\end{figure}
Note that $\dot{x_4}=x_1\dot{x}_2$ and  $\dot{x_5}=x_1\dot{x}_3$ which gives $x_4=\iint_{D_1}dx_1\wedge dx_2$ and $x_5=\iint_{D_2}dx_1\wedge dx_3$ (see Figure~\ref{fig:ds}) and proves the following.
\begin{corollary} 
Geodesics of the sub-Riemannian structure $(\TT, \DD^3, g)$ are horizontal lifts of helix arcs, with straight lines appearing as a limiting case. The coordinates $x_4$ and $x_5$ measure the signed areas associated with the projections of these geodesics onto the $x_1x_2$- and $x_1x_3$-planes, respectively.
\end{corollary} 

\subsubsection{Distribution $\DD^4$}

On the distribution $\DD^4$, the sub-Riemannian metric is given at the identity by $g^{ij}(e) = \diag(1,1,1,1,0).$ 
Using formulas~\eqref{coef:HamiltonsEq}, the metric $g^{ij}$ in coordinates takes the form
\begin{align*}
g^{ij}(p) = \begin{bmatrix}
1 & 0 & 0 & 0 & 0 \\
0 & 1 & 0 & x_1 & 0 \\
0 & 0 & 1 & 0 & x_1 \\
0 & x_1 & 0 & 1+ x_1^2 & 0 \\
0 & 0 & x_1 & 0 & x_1^2
\end{bmatrix}.
\end{align*}

This sub-Riemannian metric is compatible with the standard left-invariant Riemannian metric $\langle\cdot,\cdot\rangle$ on $\TT$ which in coordinates is represented by a symmetric matrix $S$ from~\eqref{eq:S}.

Calculation in this case gives us the Hamiltonian
\begin{align*}
H(\mathbf{x},\mathbf{\lambda} )=\frac{1}{2} \left(\lambda _1^2+(\lambda _2+\lambda_4 x_1)^2+\left(\lambda
   _3+\lambda _5 x_1\right)^2+\lambda _4^2\right),
\end{align*}

and Hamilton system:
\begin{align} \label{system:Ham4D}
  \dot{x}_1 &= \lambda_1, &  \dot{\lambda} _1 &= -\lambda_2\lambda_4-\lambda_3\lambda_5-x_1(\lambda_4^2+\lambda_5^2),\notag\\
  \dot{x}_2 &= \lambda_2+x_1\lambda_4, &  \dot{\lambda} _2 &= 0,\notag\\
  \dot{x}_3 &= \lambda_3+x_1\lambda_5, & \dot{\lambda} _3 &= 0,\\
  \dot{x}_4 &= x_1\lambda_2+(1+x_1 ^2)\lambda_4, &   \dot{\lambda} _4 &= 0,\notag\\
  \dot{x}_5 &= x_1\lambda_3+x_1^2\lambda_5, &  \dot{\lambda} _5 &= 0.\notag
\end{align}
Immediately we notice that this Hamiltonian system is completely integrable and the first integrals are given by~\eqref{eq:integrali}.
\begin{theorem}
    Hamiltonian systems~\eqref{system:Ham3D} and~\eqref{system:Ham4D} have the same first integrals given by~\eqref{eq:integrali}.
\end{theorem}

We also observe that systems~\eqref{system:Ham4D} and~\eqref{system:Ham3D} differ only in the equation for $\dot{x}_4$. Hence, following the same procedure as for the distribution $\DD^3$, we obtain the coordinate functions of the Hamiltonian characteristic curves. These coincide with those stated in Theorem~\ref{geodesics:D3}, except for the fourth coordinate, which is given by:
\begin{align*}
x_4(t)&= \tilde{x}_4(t)+t\omega\sin{\phi},
\end{align*}
where $\tilde{x}_4(t)$ denotes the fourth coordinate of the sub-Riemannian geodesic $\gamma(t)$ obtained in Theorem~\ref{geodesics:D3}.

Theorem~\ref{thm: helices_D3} extends to the case of a 4-dimensional distribution, with the proof remaining unchanged. Consequently, the following corollary is obtained.
\begin{corollary}
 Geodesics of the sub-Riemannian structure $(\TT,\DD^4,g)$ are the horizontal lifts of helix arcs, with straight lines arising as a limiting case. In particular, they are generalized helices.
\end{corollary}

The LL systems corresponding to the distributions $\DD^3$ and $\DD^4$ are highly similar. In contrast, as we shall see in the following sections, the corresponding LR systems differ substantially.

\subsection{LR systems}

The group $\TT$ also acts on itself by right translations $R_g:\TT\rightarrow \TT, \  R_g(h)=hg$ for $g,h \in \TT$. For $g=( x_1, x_2, x_3, x_4, x_5)\in \TT$, this action is given by:
\begin{align*}
dR_g(e_1)=e_1+x_2 e_4+x_3e_5, \ dR_g(e_2)=e_2, \ dR_g(e_1)=e_3, \ dR_g(e_4)=e_4, \ dR_g(e_5)=e_5. 
\end{align*}
The right invariant distributions generated by subspaces $\mathcal{L}_3$ and $\mathcal{L}_4$ are, respectively, 3-dimensional $\DD^3_R=dR_g \mathcal{L}_3$ and 4-dimensional $\DD^4_R=dR_g \mathcal{L}_4$. 

\subsubsection{Distribution $\DD^3_R$}

Sub-Riemannian metric on $\DD^3_R$:
\begin{align}\label{eq:LRmetrika3}
    g^{ij}(p)=\frac{1}{1{+}x_1^2{+}x_2^2{+}x_3^2}\begin{bmatrix}
        1 {+} x_1^2 & x_1 x_2 & x_1 x_3& (1 {+} x_1^2) x_2 & (1 {+} x_1^2) x_3\\
        x_1 x_2 & 1 {+} x_2^2 {+} \frac{x_3^2}{1 {+} x_1^2} & \frac{x_1^2 x_2 x_3}{1 {+} x_1^2} & x_1 x_2^2 & x_1 x_2 x_3\\
 x_1 x_3& \frac{x_1^2 x_2 x_3}{1 {+} x_1^2}& 1 {+} \frac{x_2^2}{1 {+} x_1^2} {+} x_3^2& 
  x_1 x_2 x_3& x_1 x_3^2\\
 (1 {+} x_1^2) x_2& x_1 x_2^2& x_1 x_2 x_3& (1 {+} x_1^2) x_2^2& (1 {+} x_1^2) x_2 x_3\\
 (1 {+} x_1^2) x_3& x_1 x_2 x_3& x_1 x_3^2& (1 {+} x_1^2) x_2 x_3& (1 {+} x_1^2) x_3^2
    \end{bmatrix}
\end{align}

The corresponding Hamiltonian takes the form:
\begin{align}\label{eq:ham3xl}
    H(\mathbf{x},\mathbf{\lambda} )=\frac{\left((1+x_1^2)(\lambda_1+x_2\lambda_4+x_3\lambda_5)+x_1(x_2\lambda_2+x_3\lambda_3)\right)^2}{2(1+x_1^2)(1{+}x_1^2{+}x_2^2{+}x_3^2)}+\frac{\lambda_2^2+\lambda_3^2}{2(1+x_1^2)}.
\end{align}

It is evident that the coordinates $x_4$ and $x_5$ are cyclic and that the Hamiltonian system has three obvious first integrals: $I_0=H$, $I_1=\lambda_4$ and $I_2=\lambda_5$. We can reduce the flow to the level set $\lambda_4=C_1$, $\lambda_5=C_2$.
Hence by setting
\begin{align*}
    \nu_1= \lambda_1+x_2 \lambda_4+x_3 \lambda_5,\quad \nu_2=\lambda_2,\quad \nu_3=\lambda_3,
\end{align*}
 we reduced the problem to the Hamiltonian system on the 6-dimensional symplectic manifold with the Poisson structure:
\begin{align*}
    \{x_1, \nu_1\}&=\{x_2, \nu_2\}=\{x_3, \nu_3\}=1,\\
    \{x_1, \nu_2\}&=\{x_1, \nu_3\}=\{x_2, \nu_1\}=\{x_2, \nu_3\}=\{x_3, \nu_1\}=\{x_3, \nu_2\}=0,\\
    \{\nu_1, \nu_2\}&=C_1,\quad\{\nu_1, \nu_3\}=C_2,\quad \{\nu_2, \nu_3\}=0.
\end{align*}
For $C_1^2+C_2^2\neq 0$, the restriction of this system corresponds to the motion of a charged particle in a constant magnetic field. 

The Hamiltonian~\eqref{eq:ham3xl} takes the form:
\begin{align}\label{eq:ham3xlrest}
     H_C(\mathbf{x},\mathbf{\nu} )=\frac{1}{2(1+x_1^2)}\left(\frac{\left((1+x_1^2)\nu_1+x_1(x_2\nu_2+x_3\nu_3)\right)^2}{1{+}x_1^2{+}x_2^2{+}x_3^2}+\nu_2^2+\nu_3^2
     \right).
\end{align}
The Hamilton's system is:
\begin{align}\label{LR:Hamsyst3}
    \begin{split}
        \dot{x}_1&=\{x_1,H_C\}=W,\\
        \dot{x}_2&=\{x_2,H_C\}=\frac{\nu_2+x_1x_2 W}{1+x_1^2},\\
        \dot{x}_3&=\{x_3,H_C\}=\frac{\nu_3+x_1x_3 W}{1+x_1^2},\\
        \dot{\nu}_1&=\{\nu_1,H_C\}=\frac{-x_1 (x_2^2 + x_3^2) W^2 -(1-x_1^2) (x_2 \nu_2 + x_3 \nu_3) W +x_1 (\nu_2^2 + \nu_3^2)}{(1 + x_1^2)^2}\\
        &\phantom{=\{\nu_1,H_C\}=}+\frac{C_1(x_1 x_2 W + \nu_2)+C_2(x_1 x_3 W + \nu_3)}{1 + x_1^2},\\
        \dot{\nu}_2&=\{\nu_2,H_C\}=-W\left(C_1+\frac{x_1\nu_2-x_2 W}{1+x_1^2}\right),\\
        \dot{\nu}_3&=\{\nu_3,H_C\}=-W\left(C_2+\frac{x_1\nu_3-x_3 W}{1+x_1^2}\right),
    \end{split}
\end{align}
where $W=\dfrac{(1+x_1^2)\nu_1+x_1(x_2\nu_2+x_3\nu_3)}{1+x_1^2+x_2^2+x_3^2}$.

Combining the above considerations yields the following proposition.
\begin{proposition}\label{prop:magnetic1}
    The reduced Hamiltonian system~\eqref{LR:Hamsyst3} with the Hamiltonian function~\eqref{eq:ham3xlrest}, obtained by symplectic reduction of the LR system at the momentum level $\lambda_4=C_1$, $\lambda_5=C_2$ is equivalent to the Hamiltonian system describing the motion of a charged particle in the constant magnetic field $F=C_1\, dx_1\wedge dx_2+C_2\, dx_1\wedge dx_3$.
\end{proposition}

The previous system can be further simplified by introducing the following coordinates $(r,\theta, \varphi)$:
\begin{align*}
    x_1=\tan\theta,\quad x_2=r\frac{\cos\varphi}{\cos\theta},\quad x_3=r\frac{\sin\varphi}{\cos\theta},
\end{align*}
with $r\geq 0$, $\theta\in\left(-\frac{\pi}{2},\frac{\pi}{2}\right)$, $\varphi\in[0,2\pi)$, and its associated coordinates $(p_r, p_\theta, p_\varphi)$:
\begin{align*}
    p_r=\frac{\sqrt{1+x_1^2}}{\sqrt{x_2^2+x_3^2}}(x_2\nu_2+x_3\nu_3)\quad p_\theta=\nu_1+x_1(x_1\nu_1+x_2\nu_2+x_3\nu_3),\quad p_\varphi=-x_3 \nu_2+x_2 \nu_3.
\end{align*}
Now, introduce the shifted angle $\psi$ as:
\begin{align*}
    \psi=\varphi-\arctan{\frac{C_2}{C_1}},\quad p_\psi=p_\varphi.
\end{align*}
In coordinates $(\mathbf{r},\mathbf{p_r})=(r,\theta,\psi, p_r, p_\theta, p_\psi)$, the Hamiltonian function takes simpler form:
\begin{align}\label{eq:ham4pol}
    \tilde{H}_C(\mathbf{r},\mathbf{p_r})=\frac{1}{2}\left(p_r^2+\frac{p_\theta^2}{1+r^2}+\frac{p_\psi^2}{r^2}\right)\cos^4\theta,
\end{align}
with the following Poisson structure:
\begin{align}\label{eq:poiss4pol}
    \begin{aligned}
        \{r,p_r\}&=\{\theta,p_\theta\}=\{\psi,p_\psi\}=1,\\
        \{r,p_\theta\}&=\{r,p_\psi\}=\{\theta,p_r\}=\{\theta,p_\psi\}=\{\psi,p_r\}=\{\psi,p_\theta\}=0,\\
        \{p_r,p_\theta\}&=-\frac{C\cos\psi}{\cos^3\theta},\quad
         \{p_r,p_\psi\}=0,\quad
         \{p_\theta,p_\psi\}=-\frac{rC\sin\psi}{\cos^3\theta},
    \end{aligned}
\end{align}
where $C=\sqrt{C_1^2+C_2^2}$. 
Hamiltonian equations in these coordinates are given by:
\begin{align}\label{eq:Hamsys3pol}
    \begin{aligned}
        \dot{r}=\{r,\tilde{H}_c\}&=\cos^4\theta\,p_r,\\
        \dot{\theta}=\{\theta,\tilde{H}_c\}&=\frac{\cos^4\theta}{1+r^2}p_\theta,\\
        \dot{\psi}=\{\psi,\tilde{H}_c\}&=\frac{\cos^4\theta}{r^2}p_\psi,\\
        \dot{p}_r=\{p_r,\tilde{H}_c\}&=\frac{r\cos^4\theta}{(1 + r^2)^2}\,p_\theta^2+\frac{\cos^4\theta}{r^3}\,p_\psi^2-\frac{C\cos\psi\cos\theta}{1+r^2}\,p_\theta,\\
        \dot{p}_\theta=\{p_\theta,\tilde{H}_c\}&=\left(p_r^2+\frac{p_\theta^2}{1+r^2}+\frac{p_\psi^2}{r^2}\right)\cos^2\theta\sin2\theta+C\left(\cos\psi\,p_r-\frac{\sin\psi}{r}\,p_\psi\right)\cos\theta,\\
        \dot{p}_\psi=\{p_\psi,\tilde{H}_c\}&=\frac{C r\sin\psi\cos\theta}{1+r^2}\,p_\theta.
    \end{aligned}
\end{align}

If $C=0$, the system has an obvious linear first integral $I_1=p_\psi$ and on a fixed level set $I_1=c$, we obtain a system with two degrees of freedom with the Hamiltonian~\eqref{eq:ham4pol} parametrically depending on the value of this integral. If $C\neq 0$, the corresponding system admits no integrals of motion linear in the momenta $p_r$, $p_\theta$, $p_\psi$. In both cases the Hamiltonian~\eqref{eq:ham4pol} is the only quadratic integral.

In order to examine the integrability of this system, we observe the restricted dynamics on the invariant manifold $\psi=p_\psi=0$. The system~\eqref{eq:Hamsys3pol} becomes:
\begin{align}\label{eq:Hamsys3polRed}
\begin{aligned}
            \dot{r}&=\cos^4\theta\,p_r, & \dot{p}_r&=\frac{r\cos^4\theta}{(1 + r^2)^2}\,p_\theta^2-\frac{C\cos\theta}{1+r^2}\,p_\theta,\\
        \dot{\theta}&=\frac{\cos^4\theta}{1+r^2}p_\theta,&
        \dot{p}_\theta&=\left(p_r^2+\frac{p_\theta^2}{1+r^2}\right)\cos^2\theta\sin2\theta+C\cos\theta\,p_r,
\end{aligned}
\end{align}
with the reduced Hamiltonian $ \tilde{H}_{red}=\dfrac{1}{2}\left(p_r^2+\dfrac{p_\theta^2}{1+r^2}\right)\cos^4\theta$.

\begin{figure}[H]
    \centering
    \includegraphics[width=0.45\linewidth]{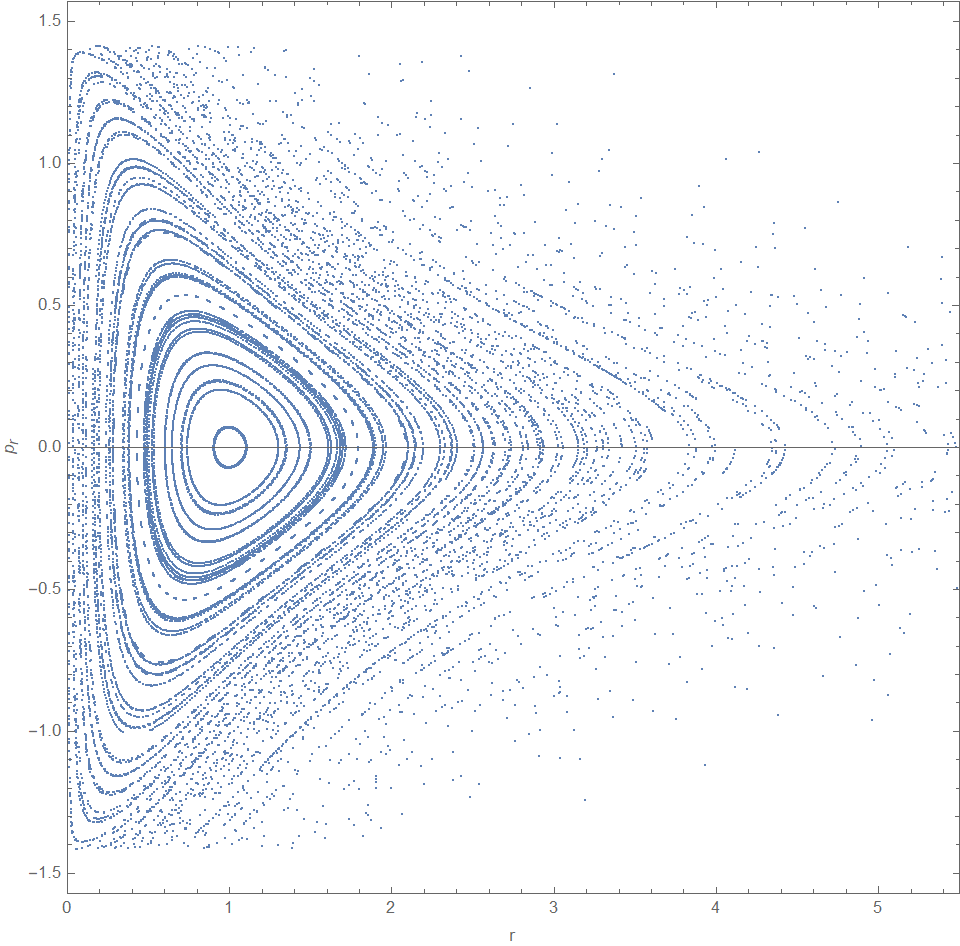}
    \caption{Poincar\'{e} section of the system~\eqref{eq:Hamsys3polRed}, with $C=3$, formed by the intersection with the plane $\theta=0$ in the coordinates $(r,p_r)$ on a level set of energy integral $\tilde{H}_{red}=1$.}
    \label{fig:poincare}
\end{figure}

Although the constructed Poincar\'{e} section exhibits regular structures when restricted to the invariant manifold (see Figure~\ref{fig:poincare}), this does not provide information about the integrability of the full system~\eqref{eq:Hamsys3pol}. In fact, numerical experiments indicate that generic trajectories are not confined to a compact region of phase space and exhibit escape behaviour in the radial variable. As a consequence, recurrent dynamics and well-defined Poincar\'{e} intersections occur only for special initial conditions, in particular those lying (or remaining close) to the invariant subsystem~\eqref{eq:Hamsys3polRed}. Hence, the observed quasi-periodic behaviour is a property of this restricted subsystem rather than evidence for integrability of the full system. 

For the purposes of the Morales–Ramis analysis, in the sequel, we regard the phase-space variables as complex-valued.

\begin{theorem}\label{thm:LR3}
    The Hamiltonian system~\eqref{eq:Hamsys3pol} with Hamiltonian~\eqref{eq:ham4pol} is not Liouville integrable. Consequently, the geodesic flow of the sub-Riemannian metric~\eqref{eq:LRmetrika3} on the right-invariant distribution $\DD^3_R$ is a Hamiltonian system with the Hamiltonian function~\eqref{eq:ham3xl} that is not integrable in Liouville sense.
\end{theorem} 
\begin{proof}
  Consider the invariant manifold $\mathcal{M}=\{\theta=p_\theta=p_\psi=0,\ \psi=\frac{\pi}{2}\}$ of the system~\eqref{eq:Hamsys3pol}. The restricted Hamiltonian system on $\mathcal{M}$ is simply $\dot{r}=p_r,\ \dot{p}_r=0$, hence the solution is 
   \begin{align}\label{eq:rr}
       r(t)=at+b,\quad a\neq 0.
   \end{align}
   
   Next, we linearize the Hamiltonian system along this particular solution by introducing infinitesimal perturbations $(\theta, \psi)$ in the directions transverse to the corresponding phase curve:
   \begin{align*}
       \dot\theta=\frac{p_\theta}{1+r^2}, \quad \dot\psi=\frac{p_\psi}{r^2},\quad \dot p_\theta=2a^2\theta-C\left(a\psi+\frac{p_\psi}{r}\right),\quad \dot p_\psi=\frac{Cr}{1+r^2}p_\theta.
   \end{align*}
This yields the normal variational equations:
\begin{align*}
   (1+r^2)\ddot\theta+2a r\dot\theta-2a^2\theta=-C(r\dot\psi+a\psi),\quad
    r^2\ddot\psi+2ar\dot\psi=Cr\dot\theta,
\end{align*}
with $r=r(t)$ being the solution~\eqref{eq:rr}. We can further simplify this by putting $\Theta(r)=\theta(t)$ and $\Psi(r)=\psi(t)$:
\begin{align}\label{eq:nve2}
    (1+r^2)\ddot\Theta+2 r\dot\Theta-2\Theta=-\frac{C}{a}(r\dot\Psi+\Psi),\quad
    r\ddot\Psi+2\dot\Psi=\frac{C}{a} \dot\Theta.
\end{align}
Now, we distinguish between two cases.
\begin{enumerate}[wide, topsep=0pt, itemsep=1pt, labelindent=0pt, label={Case \arabic*:}]
    \item\label{case:c0} $C=0$. The system~\eqref{eq:nve2} naturally decouples and can be solved explicitly:
    \begin{align*}
        \Theta=c_1 r+c_2(1+r\arctan r),\quad \Psi = c_3+\frac{c_4}{r}, \quad c_1,c_2,c_3,c_4\in\RR.
    \end{align*}
    This means that the corresponding identity component of the differential Galois group is abelian. This is necessary condition for integrability of the system, but not sufficient.  
    
    \item\label{case:cc} $C\neq0$. From the second equation in~\eqref{eq:nve2} we get that $\Theta=\frac{a}{C}(r\dot\Psi+\Psi)$,
    hence the first equation becomes:
    \begin{align}\label{eq:nve3}
  (1+r^2)\ddot\Theta+2 r\dot\Theta+(d-2)\Theta=0,
    \end{align}
where $d=\frac{C^2}{a^2}>0$. Now, put $r=i x$ and $y(x)=\Theta(r)$. The equation~\eqref{eq:nve3} reduces to Legendre differential equation:
\begin{align*}
    (1-x^2)y''-2x y'+(2-d)y = 0.
\end{align*}
The differential Galois group of the Legendre equation is abelian only for some special values of parameter $d$ (see~\cite{K69,S1872}). Hence, for generic $d$, $G=G^\circ=\SL$.
\end{enumerate}

In view of Theorem~\ref{thm:MR}, the preceding results imply that the Hamiltonian system~\eqref{eq:Hamsys3pol} is not Liouville integrable in the class of meromorphic first integrals.

Although the restriction of the dynamics to certain special levels may exhibit integrable behaviour, this phenomenon is exceptional and does not extend to the full system. Thus, the Hamiltonian system is generically non-integrable.
\end{proof}

\subsubsection{Distribution $\DD^4_R$}

Sub-Riemannian metric on $\DD^4_R$ is given by the following:
\begin{align*}
    g^{ij}(p)=\begin{bmatrix}
       \frac{1+x_1^2}{1+x_1^2+x_3^2} & 0 & \frac{x_1 x_3}{1+x_1^2+x_3^2} & 0 & \frac{\left(1+x_1^2\right) x_3}{1+x_1^2+x_3^2} \\
 0 & 1 & 0 & x_1 & 0 \\
 \frac{x_1 x_3}{1+x_1^2+x_3^2} & 0 & \frac{1+x_3^2}{1+x_1^2+x_3^2} & 0 & \frac{x_1 x_3^2}{1+x_1^2+x_3^2} \\
 0 & x_1 & 0 & 1+x_1^2 & 0 \\
 \frac{\left(1+x_1^2\right) x_3}{1+x_1^2+x_3^2} & 0 & \frac{x_1 x_3^2}{1+x_1^2+x_3^2} & 0 & \frac{\left(1+x_1^2\right) x_3^2}{1+x_1^2+x_3^2} 
    \end{bmatrix}
\end{align*}

The corresponding Hamiltonian takes the form:
\begin{align}\label{eq:ham4xl}
    H(\mathbf{x},\mathbf{\lambda} )=\frac{(x_1 (\lambda_1+x_3 \lambda_5)+x_3\lambda_3)^2+\lambda_3^2+(\lambda_1+x_3\lambda_5)^2}{2(1+x_1^2+x_3^2)}+\frac{(\lambda_2+x_1\lambda_4)^2+\lambda_4^2}{2},
\end{align}
with Hamiltonian equations:
\begin{align}\label{eq:ham4xleqs}
\begin{aligned}
    \dot{x}_1&=\{x_1, H\}=\frac{(\lambda_1+x_3\lambda_5)+x_1(x_1\lambda_1+x_3(\lambda_3+x_1\lambda_5))}{1+x_1^2+x_3^2},\\
    \dot{x}_2&=\{x_2, H\}=\lambda_2+x_1\lambda_4,\\
    \dot{x}_3&=\{x_3, H\}=\frac{\lambda_3+x_3(x_3\lambda_3+x_1(\lambda_1+x_3\lambda_5))}{1+x_1^2+x_3^2},\\
    \dot{x}_4&=\{x_4, H\}=\lambda_4+x_1(\lambda_2+x_1\lambda_4),\\
    \dot{x}_5&=\{x_5, H\}=\frac{x_3((\lambda_1+x_3\lambda_5)+x_1(x_1\lambda_1+x_3(\lambda_3+x_1\lambda_5)))}{1+x_1^2+x_3^2},\\
     \dot{\lambda}_1&=\{\lambda_1, H\}=-\frac{x_1(x_1\lambda_3-x_3(\lambda_1+x_3\lambda_5))^2}{(1+x_1^2+x_3^2)^2}+\frac{\lambda_3(x_1\lambda_3-x_3(\lambda_1+x_3\lambda_5))}{1+x_1^2+x_3^2}-(\lambda_2+x_1\lambda_4)\lambda_4,\\
     \dot{\lambda}_2&=\{\lambda_2, H\}=0,\\
     \dot{\lambda}_3&=\{\lambda_3, H\}=-\frac{(\lambda_5-x_3\lambda_1+x_1(\lambda_3+x_1\lambda_5))((\lambda_1+x_3\lambda_5)+x_1(x_1\lambda_1+x_3(\lambda_3+x_1\lambda_5)))}{(1+x_1^2+x_3^2)^2},\\
     \dot{\lambda}_4&=\{\lambda_4, H\}=0,\\
     \dot{\lambda}_5&=\{\lambda_5, H\}=0.
\end{aligned}
\end{align}
Obviously, we have three cyclic coordinates: $x_2$, $x_4$ and $x_5$. 
Hence, by setting:
\begin{align*}
    \nu_1=\lambda_1+ x_3 \lambda_5,\quad \nu_k=\lambda_k,\quad k=2,\ldots,5,
\end{align*}
 we can reduce the flow to the level set $\lambda_2=C_1$, $\lambda_4=C_2$, $\lambda_5=C_3$ and project this restriction onto the plane $(x_1, x_3, \nu_1, \nu_3)$.
Now, the problem is reduced to the system with Hamiltonian:
\begin{align}\label{eq:ham4xn}
    H_C(x_1,x_3,\nu_1,\nu_3)=\frac{1}{2}\left(\nu_1^2+\nu_3^2-\frac{(x_1\nu_3-x_3\nu_1)^2}{1+x_1^2+x_3^2}+(C_1+C_2x_1)^2\right).
\end{align}
Note that the original Hamiltonian~\eqref{eq:ham4xl} contains the additive term $\lambda_4^2$. In the reduced Hamiltonian, this term contributes only an overall constant shift and therefore has no effect on the dynamics; consequently, it is omitted in~\eqref{eq:ham4xn}. 
The Poisson structure is given by:
\begin{align}
    \{x_1, x_3\}=0,\quad \{x_1,\nu_1\}=\{x_3, \nu_3\}=1,\quad \{x_1,\nu_3\}=\{x_3,\nu_1\}=0,\quad \{\nu_1, \nu_3\}=C_3.
\end{align}
The system of Hamilton's equations:
\begin{align}\label{LR:Hamsyst4}
    \begin{split}
        \dot{x}_1&=\{x_1,H_C\}=\frac{\nu_1+x_1(x_1\nu_1+x_3\nu_3)}{1+x_1^2+x_3^2},\\
        \dot{x}_3&=\{x_3,H_C\}=\frac{\nu_3+x_3(x_1\nu_1+x_3\nu_3)}{1+x_1^2+x_3^2},\\
        \dot{\nu}_1&=\{\nu_1,H_C\}=\frac{(\nu_3+x_3(x_1\nu_1+x_3\nu_3))(x_1\nu_3-x_3\nu_1+C_3(1+x_1^2+x_3^3))}{(1+x_1^2+x_3^3)^2}-C_2(C_1+C_2 x_1),\\
        \dot{\nu}_3&=\{\nu_3,H_C\}=-\frac{(\nu_1+x_1(x_1\nu_1+x_3\nu_3))(x_1\nu_3-x_3\nu_1+C_3(1+x_1^2+x_3^3))}{(1+x_1^2+x_3^3)^2}.
    \end{split}
\end{align}

The above considerations yield the following proposition.
\begin{proposition}\label{prop:magnetic2}
    The reduced Hamiltonian system~\eqref{LR:Hamsyst4} is equivalent to the Hamiltonian system describing the motion of a charged particle in the constant magnetic field $F=C_3\, dx_1\wedge dx_3$.
\end{proposition}

\begin{figure}[H]
    \centering
    \includegraphics[width=0.45\linewidth]{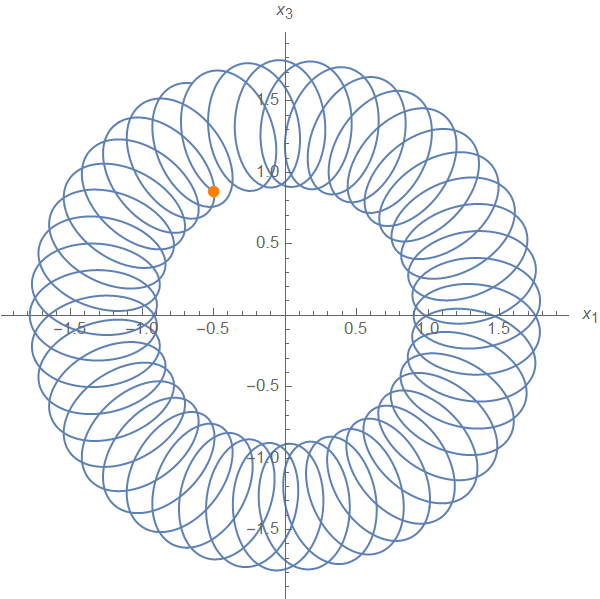}
    \includegraphics[width=0.45\linewidth]{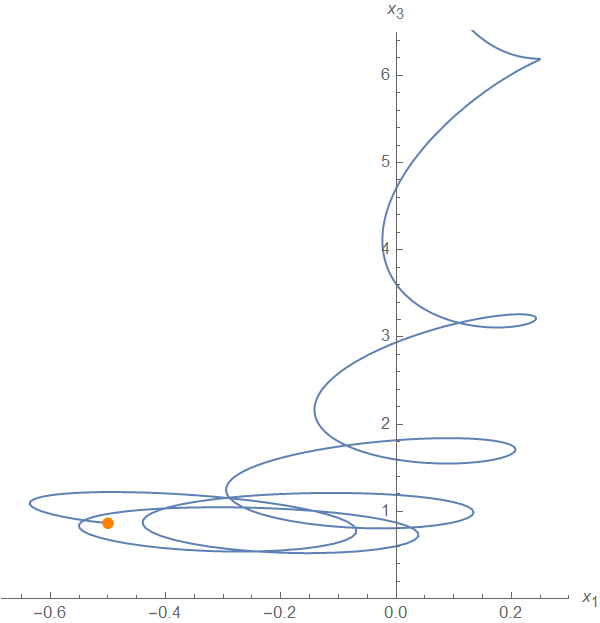}
    \caption{Geodesic lines of the system~\eqref{LR:Hamsyst4} on a level set of energy integral $H=1$, with $\lambda_2=\frac{1}{\sqrt 2}$,  $\lambda_5=5$ and (a) $\lambda_4=0$ (\textit{left}) and (b) $\lambda_4={2\sqrt 2}$ (\textit{right}), with the same initial conditions $x_1(0)=-\frac{1}{2}$, $x_3(0)=\frac{\sqrt 3}{2}$.}
    \label{fig:geod}
\end{figure}

Figure~\ref{fig:geod} illustrates geodesic trajectories of the restricted system~\eqref{LR:Hamsyst4} originating from the same initial point, contrasting the two distinguished cases: $C_2=0$ (\textit{left}) and $C_2\neq 0$ (\textit{right}). 

The figure provides strong evidence that the  system~\eqref{LR:Hamsyst4} may be integrable in some cases and non-integrable in others. However, we need more rigorous proof of this observation.

\begin{proposition}\label{prop:C2int}
    The system obtained by symplectic reduction of the LR system~\eqref{eq:ham4xleqs} at the momentum level $\lambda_4=0$ is equivalent to the Hamiltonian system with the Hamiltonian function:
    \begin{align*}
            H_0(\mathbf{x},\mathbf{\lambda} )=\frac{(x_1 (\lambda_1+x_3 \lambda_5)+x_3\lambda_3)^2+\lambda_3^2+(\lambda_1+x_3\lambda_5)^2}{2(1+x_1^2+x_3^2)}+\frac{\lambda_2^2}{2},
    \end{align*}
integrable via first integrals:
\begin{align*}
    I_1=\lambda_2,\quad I_2=\lambda_5,\quad  I_3=-x_3\lambda_1+ x_1\lambda_3+\frac{1}{2}(x_1^2-x_3^2)\lambda_5,\quad I_4=H_0,
\end{align*}
which are commutative and functionally independent almost everywhere.
\end{proposition}
\begin{proof}
As previously noticed, $I_1=\lambda_2$ and $I_2=\lambda_5$ are the first integrals. Hence, we can further reduce this system to the level set $\lambda_2=C_1$, $\lambda_5=C_3$ and introduce the polar coordinates
$(r,\theta,p_r,p_\theta)$:
\begin{align*}
    r=\sqrt{x_1^2+x_3^2},\quad\theta=\arctan{\frac{x_3}{x_1}},\quad p_r =\frac{x_1 \nu_1 + x_3 \nu_3}{\sqrt{x_1^2 + x_3^2}},\quad p_\theta= -x_3\nu_1 + x_1\nu_3,
\end{align*}
The Hamiltonian
\begin{align*}   
\tilde{H}_C(r,\theta,p_r,p_\theta)=\frac{1}{2}\left(p_r^2+\frac{p_\theta^2}{r^2(1+r^2)}\right)+\frac{C_1^2}{2},
\end{align*}
consists only of kinetic part and the system:
\begin{align*}
        \dot{r}&=\{r,\tilde{H}_C\}=p_r,&  \dot{p}_r&=\{p_r,\tilde{H}_C\}=\frac{(1+2r^2)}{r^3(1+r^2)^2}p_\theta^2+\frac{C_3}{r(1+r^2)}p_\theta,\\
        \dot{\theta}&=\{\theta,\tilde{H}_C\}=\frac{p_\theta}{r^2(1+r^2)},&
        \dot{p}_\theta&=\{p_\theta,\tilde{H}_C\}=-C_3r p_r.
\end{align*}
has two obvious integrals:
\begin{align}\label{eq:heispol}
    I_3=p_\theta+\frac{C_3}{2}r^2\quad\text{and}\quad I_4=p_r^2+\frac{p_\theta^2}{r^2(1+r^2)}.
\end{align}
In the original coordinates, the integral $I_3$ becomes
\begin{align*}
    I_3=-x_3\lambda_1+ x_1\lambda_3+\frac{\lambda_5}{2}(x_1^2-x_3^2),
\end{align*}
while $I_4$ is the Hamiltonian itself.
These integrals are involutive and functionally independent almost everywhere. 
\end{proof}

Observe that any system reduced to the level set $\lambda_2=C_1$, $\lambda_4=0$ exhibits the same properties as the system restricted to $\lambda_2=\lambda_4=0$, since the constant $C_1$ enters the Hamiltonian only through an additive term. Furthermore, this reduced system possesses a remarkable feature, as stated in the following proposition.
\begin{proposition}\label{prop:heis}
    The system obtained by symplectic reduction of the LR system~\eqref{eq:ham4xleqs} at the momentum level $\lambda_2=\lambda_4=0$ is equivalent to the LR Hamiltonian system on the Heisenberg group $H_3$. 
\end{proposition}

\begin{proof}
    In~\cite{T97} was shown that the geodesic flow of the sub-Riemannian metric on $H_3$ corresponding to the left-invariant Riemannian metric $ds^2=dx^2+(1+x^2)dy^2+dz^2-2x\,dx\,dz$ and the right-invariant distribution $R_{g*}\DD=\RR\langle\frac{\partial}{\partial x}+y\frac{\partial}{\partial z}, \frac{\partial}{\partial y}\rangle$ constitutes the Hamiltonian system with Hamiltonian function:
    \begin{align*}
        H(\mathbf{x},\mathbf{\mu})&=\frac{(x(\mu_1+y\mu_3)+y\mu_2)^2+\mu_2^2+(\mu_1+y\mu_3)^2}{2(1+x^2+y^2)},
    \end{align*}
    which is exactly the Hamiltonian~\eqref{eq:ham4xl} with $\lambda_2=\lambda_4=0$.
\end{proof}

For the non-constant potential, we can consider the simplest case with $C_1=C_3=0$ and $C_2^2=\epsilon>0$. Now, we make linear canonical transformation:
\begin{align*}
   q_1=\frac{1}{\sqrt 2}(x_3-i x_1), \quad  q_2=\frac{1}{\sqrt 2}(x_3+i x_1), \quad p_1=\frac{1}{\sqrt 2}(\nu_3+i \nu_1),\quad
   p_2=\frac{1}{\sqrt 2}(\nu_3-i \nu_1).
\end{align*}
In the new variable Hamiltonian~\eqref{eq:ham4xn} reads:
\begin{align*}
    H_\epsilon^\mathbb{C}(q_1,q_2,p_1,p_2)=\frac{1}{2}\left(\frac{2 p_1 p_2+(q_1 p_1+ q_2 p_2)^2}{1+2 q_1 q_2}-\frac{\epsilon}{2}(q_1-q_2)^2\right),
\end{align*}
and the equations of motion are
\begin{align}\label{eq:complexeq}
\begin{split}
   \dot{q}_1&=\frac{q_1^2 p_1+(1+q_1 q_2)p_2}{1+2 q_1 q_2}, \quad \dot{p}_1=-\frac{(q_1 p_1-q_2 p_2)(q_2^2 p_2+(1+q_1 q_2)p_1)}{(1+2 q_1 q_2)^2}+\frac{\epsilon}{2}(q_1-q_2),\\
   \dot{q}_2&=\frac{q_2^2 p_2+(1+q_1 q_2)p_1}{1+2 q_1 q_2}, \quad \dot{p}_2=\phantom{-}\frac{(q_1 p_1-q_2 p_2)(q_1^2 p_1+(1+q_1 q_2)p_2)}{(1+2 q_1 q_2)^2}-\frac{\epsilon}{2}(q_1-q_2).  
\end{split}
\end{align}
Fix the invariant manifold $q_1+q_2=0$, $p_1+p_2=0$ and consider the reduced dynamics on it. Set $q_1=-q_2=q$ and $p_1=-p_2=p$. The Hamiltonian restricted to this plane becomes
\begin{align*}
    H_{red}^\mathbb{C}=-p^2-\epsilon q^2
\end{align*}
and the particular solution of the corresponding system~\eqref{eq:complexeq} can be written explicitly:
\begin{align}\label{eq:part}
    q(t)=a \cos(t\sqrt{\epsilon}),\quad p(t)=a\sqrt{\epsilon}\sin(t\sqrt{\epsilon}),\quad a\neq 0.
\end{align}
Now, put
   $ q_1=q+\eta,\quad q_2=-q+\eta,\quad p_1=p+\pi,\quad p_2=-p+\pi$.
The quadratic transverse Hamiltonian becomes
\begin{align*}
    H_{trans}^\mathbb{C}=\frac{\pi^2+4 q p\eta\pi+2p^2\eta^2}{1 - 2 q^2}.
\end{align*}
Since $p_1\, dq_1+ p_2\, dq_2=2 p\, dq + 2\pi\, d\eta$, the normal variables have symplectic form $2d\pi\wedge d\eta$. Therefore
\begin{align*}
    \dot\eta=\frac{1}{2}\frac{\partial H_{trans}^\mathbb{C}}{\partial \pi}=\frac{\pi+2pq\eta}{1- 2 q^2},\qquad
    \dot\pi=-\frac{1}{2}\frac{\partial H_{trans}^\mathbb{C}}{\partial \eta}=\frac{-2p(q\pi+p\eta)}{1- 2 q^2}.
\end{align*}
Using $\dot q=-p$, $\dot p=\epsilon q$ and $p^2+\epsilon q^2=a^2 \epsilon$, this can be combined to a second-order equation:
\begin{align}\label{eq:nve}
    \ddot\eta(t) + \frac{4q(t) p(t)}{1- 2 q^2(t)}\dot\eta(t)+\frac{4a^2 \epsilon-6\epsilon q^2(t)}{1- 2 q^2(t)}\eta(t)=0,
\end{align}
where $q(t)$ and $p(t)$ are given by~\eqref{eq:part}. Now, set $\theta(s)=\eta(t)\sqrt{1-2q^2(t)}$, where $s=t\sqrt\epsilon$. The normal variational equation~\eqref{eq:nve} becomes
\begin{align}\label{eq:hill}
\ddot\theta(s)=-\frac{a^2(1 + 2 a^2 + (-5 + 4 a^2) \cos2s + 2 a^2 \cos4 s)}{(1 - 2 a^2\cos^2 s)^2}\theta(s)=R(s)\theta(s).
\end{align}
Hence, the equation~\eqref{eq:hill} is Hill equation defined by the $\pi$-periodic rational function $R(s)$. Since $R$ is basically a function of $\cos^2 s$, we can make a further adjustment by setting  $\theta(s)=z^{-1/4}(1-z)^{-1/4}\phi(z)$ with $z=\cos^2 s$:
\begin{align*}
\ddot\phi(z)=\frac{-3 + 4(1 - 3 a^2) z + 4(-1 + 12 a^2 + 5 a^4) z^2 - 
 8 a^2 (3 + 10 a^2) z^3 + 48 a^4 z^4}{16 z^2 (1 - z)^2 (1 -2 a^2 z)^2}\phi(z).
\end{align*}
This equation has exactly four regular singular points: $z=0$, $z=1$, $z=\frac{1}{2a^2}$ and $z=\infty$. Kovacic algorithm~\cite{K86} finds no Liouvillian solution for generic $a\neq 0$, hence the differential Galois group is $\SL$. As $\SL$ is connected, its identity component is the group itself, which is non-abelian.  

The previous analysis made us possible to formulate the following result.
\begin{theorem}\label{th:LR4}
Consider the reduced Hamiltonian system~\eqref{LR:Hamsyst4} with the Hamiltonian function~\eqref{eq:ham4xn}, by symplectic reduction of the LR system~\eqref{eq:ham4xleqs} at the momentum level $\lambda_2=C_1$, $\lambda_4=C_2$, $\lambda_5=C_3$. 
For $C_2=0$ the system is equivalent to the reduced Hamiltonian system on Heisenberg group $H_3$, hence it is  Liouville integrable. 
For $C_2\neq0$, this system is not  Liouville integrable.
 Consequently, the LR system~\eqref{eq:ham4xleqs} is not Liouville integrable.
\end{theorem}
\begin{proof}
If $C_2=0$, we have shown in Proposition~\ref{prop:C2int} that there exists an exceptional invariant manifold on which the dynamics simplifies. Moreover, by Proposition~\ref{prop:heis}, this system is equivalent to the Hamiltonian system on Heisenberg group $H_3$ restricted to the level set $\mu_3=const$. This system is integrable in the Liouville sense with the first integrals given by~\eqref{eq:heispol} (in polar coordinates). 

The non-integrability of the reduced system for $C_2\neq0$ follows from the Galoisian obstruction~\cite{MRR01}, since we showed that, in special case when $C_1=C_3=0$, the identity component of the Galois group of the variational equation is $\SL$, hence it is not abelian. Therefore, by Theorem~\ref{thm:MR}, the Hamiltonian system~\eqref{LR:Hamsyst4} fails to be Liouville integrable with meromorphic first integrals.

If the original Hamiltonian system were Liouville integrable, then its first integrals would restrict to the invariant level set $\lambda_2=C_1$, $\lambda_4=C_2$, $\lambda_5=C_3$ 
yielding a complete set of commuting first integrals for the reduced system. This contradicts its non-integrability. Therefore, the original Hamiltonian system is not Liouville integrable. 
\end{proof}

\section{Group $\HH$}\label{sec:H}

In this section, we consider the so-called second Heisenberg group $\HH_5$, a connected 2-step nilpotent Lie group with a one-dimensional center, whose matrix representation is given by:
\begin{align*}
\HH = \left\{
\begin{bmatrix}
1 & x_1 & x_2 & x_5 \\
0 & 1 & 0 & x_3 \\
0 & 0 & 1 & x_4 \\
0 & 0 & 0 & 1 
\end{bmatrix}: x_1, x_2, x_3, x_4, x_5 \in \mathbb{R} \right\}.
\end{align*}
The corresponding Lie algebra $\hh$, generated by $e_1,\ldots,e_5$, has the following non-trivial commutator relations:
\begin{align*}
[e_1, e_3]= [e_2, e_4]=e_5.
\end{align*}

In the sequel, we recall previously known results. While the LL systems are available in a more general setting, the LR systems were studied in detail in~\cite{PS,PAV}, albeit only for the standard metric. In this work, we extend these results to a general metric and adapt them to the 5-dimensional case.

\subsection{LL system}

In~\cite{BN13} it was proven that any left-invariant left-invariant sub-Riemannian structure $(\DD, g)$ on $\HH_5$ is isometric to exactly one of the structures $(\mathcal{H}^4, g^\sigma)$, where  $\mathcal{H}^4=\RR\langle e_1, e_2, e_3, e_4\rangle$ is the contact distribution and
\begin{align*}
 g^\sigma(e)  = D(\sigma)= \diag(\sigma,1,\sigma,1),\quad  \sigma\geq 1.
\end{align*}

The sub-Riemannian metric:
\begin{align}\label{eq:subR_l}
g^{ij}(p) = \begin{bmatrix}
 \frac{1}{\sigma} & 0 & 0 & 0 & 0 \\
 0 & 1 & 0 & 0 & 0 \\
 0 & 0 & \frac{1}{\sigma} & 0 & \dfrac{x_1}{\sigma}  \\
 0 & 0 & 0 & 1 & x_2 \\
 0 & 0 & \dfrac{x_1}{\sigma} & x_2 & \dfrac{x_1^2}{\sigma}+x_2^2
\end{bmatrix}
\end{align}
is compatible with a left-invariant Riemannian metric on $\HH_5$, whose coordinate representation is given by the symmetric matrix:

\begin{align}\label{eq:leftRim}
S = \begin{bmatrix}
    \sigma & 0 & 0 & 0 & 0\\
    0 & 1 & 0 & 0 & 0 \\
    0 & 0 & \tau x_1^2+\sigma & \tau x_1 x_2 & -\tau x_1\\
    0 & 0 & \tau x_1 x_2 & \tau x_2^2+1 & -\tau x_2\\
     0 & 0 & -\tau x_1 & -\tau x_2 & \tau
\end{bmatrix}. 
\end{align}

\begin{theorem}\label{thm:left}
\begin{enumerate}[wide, topsep=0pt, itemsep=10pt,  labelindent=0pt, label = (\alph*)]
  \item  The geodesic flow corresponding to the sub-Riemannian metric~\eqref{eq:subR_l} and the left-invariant distribution $dL_g\mathcal{H}^4$, is a Hamiltonian system on $T^*\HH_5$ with the Hamiltonian function
\begin{align*}
  H(q,\lambda)=\frac{1}{2\sigma}[\lambda_1^2+\sigma\lambda_2^2+(\lambda_{3}+x_1\lambda_{5})^2+\sigma(\lambda_{4}+x_2\lambda_{5})^2],
\end{align*} 
  integrable in a non-commutative sense via first integrals:
  \begin{gather}
  \begin{aligned}\label{eq:intLL}
I_0 &=\lambda_5,\  I_1 = \lambda_3,\ I_2 = \lambda_4,\  I_3=\lambda_1+x_3\lambda_5,\  I_4=\lambda_2+x_4\lambda_5,\\ I_5&=\lambda_1^2+(\lambda_3+x_1\lambda_5)^2,\ I_6=\lambda_2^2+(\lambda_4+x_2\lambda_5)^2,
\end{aligned}
  \end{gather}
  which are functionally independent almost everywhere, and the integrals $I_0$, $I_5$ and $I_6$ commute with all the other integrals. Moreover, the subset $\{\lambda_5\neq 0\}$ of its phase space is foliated into general helices:
  \begin{align}
    \begin{split}\label{eq:traj}
    x_1(t) &= a_1 \cos\left(\frac{\lambda_5}{\sigma} t\right) + b_1 \sin\left(\frac{\lambda_5}{\sigma} t\right) + c_1,\quad 
    x_2(t) = a_2 \cos\left(\lambda_5 t\right) + b_2 \sin\left(\lambda_5 t\right) + c_2,\\
    x_3(t) &= -b_1 \cos\left(\frac{\lambda_5}{\sigma} t\right) + a_1 \sin\left(\frac{\lambda_5}{\sigma} t\right) + d_1,\quad
    x_4=-b_2 \cos\left(\lambda_5 t\right) + a_2 \sin\left(\lambda_5 t\right) + d_2,
    \end{split}\\
    z(t)&=c_0+\frac{\left(a_1^2{+}b_1^2\right)\lambda_5}{2\sigma} t
   {+}\frac{\left(a_1^2{-}b_1^2\right)}{4}\sin
   \left(\frac{2\lambda_5}{\sigma} t\right){-}\frac{a_1 b_1}{2}  \cos \left(\frac{2\lambda_5}{\sigma} t\right){+} a_1 c_1 \sin \left(\frac{\lambda_5}{\sigma} t\right)-b_1 c_1 \cos \left(\frac{\lambda_5}{\sigma} t\right)\notag\\
   &\phantom{=c_0}+\frac{\left(a_2^2+b_2^2\right)\lambda_5}{2} t
   {+}\frac{\left(a_2^2-b_2^2\right)}{4} \sin
   \left(2\lambda_5 t\right){-}\frac{a_2 b_2}{2}  \cos \left(2\lambda_5 t\right){+} a_2 c_2 \sin \left(\lambda_5 t\right)-b_2 c_2 \cos \left(\lambda_5 t\right),\notag
\end{align}
with $a_1, a_2, b_1, b_2 ,c_0, c_1, c_2, d_1, d_2\in\RR$. 

  \item The  reduced system obtained by restricting the geodesic flow to the level set $\lambda_5=C\neq 0$ and projecting this restriction of the flow to the plane $\RR^{8}(x,y,\lambda)$ is equivalent to the Hamiltonian system describing the motion of a charged particle in the constant magnetic field $F=-C \left(\frac{1}{\sigma}\, dx_1\wedge dx_3+dx_2\wedge dx_4\right)$. This flow is superintegrable and its phase space for $C\neq0$ is foliated into trajectories~\eqref{eq:traj} that are either closed or dense on a torus $T^2$ (depending on a rational dependency of constants $C$ and $\sigma$).
    \end{enumerate}
\end{theorem}

\begin{figure}[H]
    \centering
    \includegraphics[width=0.43\linewidth]{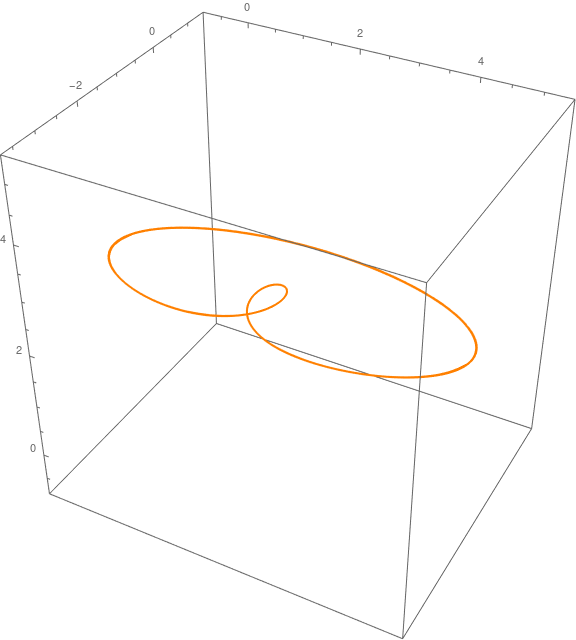}\qquad
    \includegraphics[width=0.43\linewidth]{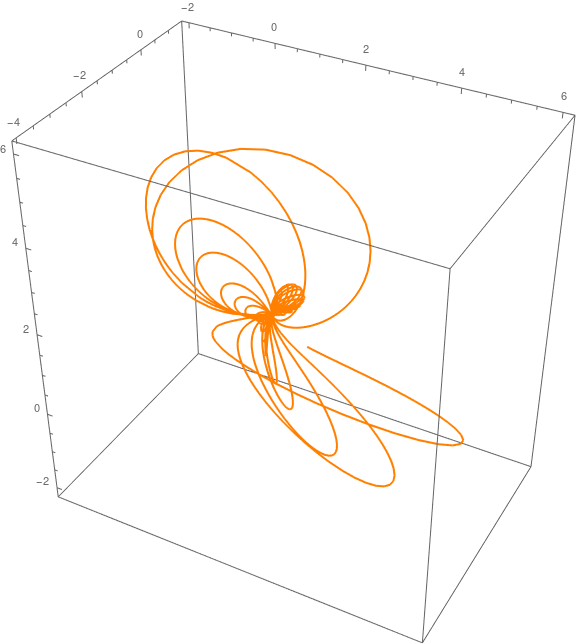}
    \caption{Stereographic projection $(x_1,x_2,x_3,x_4)\mapsto\left(\dfrac{x_1}{1-x_4},\dfrac{x_2}{1-x_4},\dfrac{x_3}{1-x_4}\right)$ of the trajectories~\eqref{eq:traj} with (a) rational dependency (\textit{left}); (b) irrational dependency (\textit{right}) of the constants $C$ and $\sigma$.}
    \label{fig:geodH5red}
\end{figure}

\subsection{LR system}

 The right-invariant distribution generated by $\mathcal{H}^4$ consists of $4$-planes $\mathcal{H}^4_R=R_{g*}\mathcal{H}^4$, where:
\begin{align*}
&R_{g*}(\frac{\partial}{\partial x_1})=\frac{\partial}{\partial x_1}+ x_3 \frac{\partial}{\partial x_5},\quad
R_{g*}(\frac{\partial}{\partial x_2})=\frac{\partial}{\partial x_2}+ x_4 \frac{\partial}{\partial x_5},\\
&R_{g*}(\frac{\partial}{\partial x_3})=\frac{\partial}{\partial x_3},\quad R_{g*}(\frac{\partial}{\partial x_4})=\frac{\partial}{\partial x_4}\quad\text{and}\quad
R_{g*}(\frac{\partial}{\partial x_5})=\frac{\partial}{\partial x_5}.
\end{align*}

The sub-Riemannian metric $g^{ij}(x_1,\ldots,x_5)$ on the right-invariant distribution $\mathcal{H}^4_R$, corresponding to the left-invariant Riemannian metric~\eqref{eq:leftRim}, is given by the matrix:
   \begin{align*}
   \setlength{\arraycolsep}{3pt}
\frac{1}{s}{\scriptstyle\begin{pmatrix}
    \dfrac{(s {-} x_3^2 \tau)}{\sigma}& {-}x_3 x_4 \tau & \dfrac{x_1 x_3 \tau}{\sigma}& 
  x_2 x_3 \tau& 
  \dfrac{x_3 (\sigma {+} x_1^2 \tau {+} 
     x_2^2 \sigma \tau)}{\sigma}\\ 
     {-}x_3 x_4 \tau & (s {-} x_4^2 \sigma \tau) & x_1 x_4 \tau & x_2 x_4 \sigma \tau & 
  x_4 (\sigma {+} x_1^2 \tau {+} x_2^2 \sigma \tau)\\ 
  \dfrac{x_1 x_3 \tau}{\sigma} & x_1 x_4 \tau & \dfrac{(s {-} x_1^2 \tau)}{\sigma} & {-}x_1 x_2 \tau & 
  \dfrac{x_1 (x_3^2 {+} x_4^2 \sigma) \tau}{\sigma}\\
  x_2 x_3 \tau & x_2 x_4 \sigma \tau & {-}x_1 x_2 \tau & (s {-} x_2^2 \sigma \tau) & 
  x_2 (x_3^2 {+} x_4^2 \sigma) \tau\\ 
 \dfrac{x_3 (\sigma {+} x_1^2 \tau {+} 
     x_2^2 \sigma \tau)}{\sigma} & 
  x_4 (\sigma {+} x_1^2 \tau {+} x_2^2 \sigma \tau) & 
  \dfrac{x_1 (x_3^2 {+} x_4^2 \sigma) \tau}{\sigma} &
  x_2 (x_3^2 {+} x_4^2 \sigma) \tau & \dfrac{(x_3^2 {+} x_4^2 \sigma) (\sigma {+} x_1^2 \tau {+} 
      x_2^2 \sigma \tau)}{\sigma}
  \end{pmatrix}},
\end{align*}
where $s=\sigma+(x_1^2+x_3^2+(x_2^2+x_4^2)\sigma)\tau$.

As a consequence, we obtain the following proposition.
\begin{proposition}
   The LR system on $T^*H_5$ is a Hamiltonian system  with the Hamiltonian function:
    \begin{gather}\label{eq:HamR}
       \begin{aligned}
           H(\mathbf{x},\mathbf{\lambda})=\frac{1}{2\sigma s}&\left(\lambda_1^2 (\sigma {+} x_1^2 \tau {+} (x_2^2 {+} x_4^2) \sigma \tau){+}\lambda_2^2 \sigma (\sigma {+} (x_1^2 {+} x_3^2 {+} x_2^2 \sigma) \tau)\right.\\
           &{+}\lambda_3^2 (\sigma {+} x_3^2 \tau {+} (x_2^2 {+} x_4^2) \sigma \tau) {+}\lambda_4^2 \sigma (\sigma {+} (x_1^2 {+} x_3^2 {+} x_4^2 \sigma) \tau)\\
           &{+}2\tau\sigma (\frac{1}{\sigma}x_1 x_3 \lambda_1\lambda_3 {+} x_2 x_3  \lambda_1\lambda_4 {+}x_1 x_4\lambda_2\lambda_3{+}x_2 x_4 \lambda_2\lambda_4 \sigma{-}x_3 x_4 \lambda_1\lambda_2{-}x_1 x_2 \lambda_3 \lambda_4)\\
           &{+}2\lambda_5((\sigma {+} (x_1^2 {+} x_2^2 \sigma) \tau)( x_3 \lambda_1 {+} x_4\lambda_2\sigma){+}(x_3^2 {+} x_4^2 \sigma)\tau(x_1\lambda_3{+}x_2\lambda_4 \sigma ))\\
           &\left.{+}\lambda_5^2 (x_3^2 {+} x_4^2 \sigma) (\sigma {+} (x_1^2 {+} 
   x_2^2 \sigma) \tau)\right),
       \end{aligned}
    \end{gather}
and canonical Poisson structure. Hamiltonian equations take the form:
\begin{gather}
\begin{aligned}\label{eq:HamReq}
    \dot{x_1} =\{x_1, H\} &= \frac{x_3 \lambda_5 \sigma {+} 
 x_3 (x_1 (x_1 \lambda_5 {+} \lambda_3) {-} x_4 \lambda_2 \sigma {+} 
    x_2 (x_2 \lambda_5 {+} \lambda_4)\sigma)\tau {+} \lambda_1 (\sigma{+} 
    x_1^2\tau {+} (x_2^2 {+} x_4^2)\sigma \tau)}{\sigma(\sigma{+}(x_1^2{+}x_3^2{+}(x_2^2{+}x_4^2)\sigma)\tau)},\\
    \dot{x_2} =\{x_2, H\} &= \frac{x_4 \lambda_5 \sigma {+} 
 x_4 (x_1 (x_1 \lambda_5 {+} \lambda_3) {-} x_3 \lambda_1 {+} 
    x_2 (x_2 \lambda_5 {+} \lambda_4)\sigma)\tau {+} \lambda_2 (\sigma{+} 
    (x_1^2 {+} x_3^2 {+} x_2^2\sigma) \tau)}{\sigma{+}(x_1^2{+}x_3^2{+}(x_2^2{+}x_4^2)\sigma)\tau},\\
    \dot{x_3} =\{x_3, H\} &= \frac{x_1 (x_3 (x_3 \lambda_5 {+} \lambda_1) {+} 
    x_4 (x_4 \lambda_5 {+} \lambda_2) \sigma {-} 
    x_2 \lambda_4 \sigma) \tau {+} \lambda_3 (\sigma {+} 
    x_3^2 \tau {+} (x_2^2 {+} x_4^2) \sigma \tau)}{\sigma(\sigma{+}(x_1^2{+}x_3^2{+}(x_2^2{+}x_4^2)\sigma)\tau)},\\
    \dot{x_4} =\{x_4, H\} &= \frac{x_2 (x_3 (x_3 \lambda_5 {+} \lambda_1) {+} 
    x_4 (x_4 \lambda_5 {+} \lambda_2) \sigma {-} 
    x_1 \lambda_3 \sigma) \tau {+} \lambda_4 (\sigma {+} 
    (x_1^2 {+} x_3^2{+} x_4^2 \sigma) \tau)}{\sigma{+}(x_1^2{+}x_3^2{+}(x_2^2{+}x_4^2)\sigma)\tau},\\
    \dot{z}=\{z, H\}&=\dot{x}_1x_3+\dot{x}_2x_4,\\
    \dot{\lambda_1}=\{\lambda_1, H\}&=-\frac{\left(x_3\lambda_1{-} x_1\lambda_3{+}( x_4\lambda_2{-} x_2\lambda_4)\sigma{+} (x_3^2{+} x_4^2\sigma)\lambda_5\right)\tau}{\sigma{+}\left(x_1^2{+} x_3^2{+} (x_2^2 {+} x_4^2)\sigma\right)\tau}\dot{x}_3,\\
    \dot{\lambda_2}=\{\lambda_1, H\}&={-}\frac{\left(x_3\lambda_1{-} x_1\lambda_3{+}( x_4\lambda_2{-} x_2\lambda_4)\sigma{+} (x_3^2{+} x_4^2\sigma)\lambda_5\right)\tau}{\sigma{+}\left(x_1^2{+} x_3^2{+} (x_2^2 {+} x_4^2)\sigma\right)\tau}\dot{x}_4,\\
    \dot{\lambda_3}=\{\lambda_3, H\}&=\left(\lambda_5{+}\frac{\left(x_3\lambda_1- x_1\lambda_3{+}( x_4\lambda_2- x_2\lambda_4)\sigma{+} (x_3^2{+} x_4^2\sigma)\lambda_5\right)\tau}{\sigma{+}\left(x_1^2{+} x_3^2{+} (x_2^2 {+} x_4^2)\sigma\right)\tau}\right)\dot{x}_1,\\
    \dot{\lambda_4}=\{\lambda_4, H\}&=\left(\lambda_5{+}\frac{\left(x_3\lambda_1- x_1\lambda_3{+}( x_4\lambda_2- x_2\lambda_4)\sigma{+} (x_3^2{+} x_4^2\sigma)\lambda_5\right)\tau}{\sigma{+}\left(x_1^2{+} x_3^2{+} (x_2^2 {+} x_4^2)\sigma\right)\tau}\right)\dot{x}_2,\\
     \dot{\lambda_5}=\{\lambda_5, H\}&=0.
\end{aligned}
\end{gather}
\end{proposition}
The explicit expression for $\dot{z}$ in~\eqref{eq:HamReq} can be obtained by direct calculation. Since it is rather lengthy and plays no further role in the analysis, its full form is omitted.

\begin{theorem}\label{thm:rightR}
   The reduced Hamiltonian system obtained by symplectic reduction of the LR system~\eqref{eq:HamReq} at the momentum level $\lambda_5=C$ is equivalent to the Hamiltonian system describing the motion of a charged particle in the constant magnetic field $F=C\left(\frac{1}{\sigma }dx_1\wedge dx_3+dx_2\wedge dx_4\right)$.  This system is completely integrable and its first integrals are given by~\eqref{eq:Hc},~\eqref{eq:lin} and~\eqref{eq:quad}.
 \end{theorem}

 \begin{proof}
From~\eqref{eq:HamReq}, it is obvious that the variable $x_5$ is cyclic. Hence, we can restrict the flow to the level set $\lambda_5 = C$ and in this way, we reduce the problem to the Hamiltonian system on the $8$-dimensional symplectic manifold. Set $q_i=x_i, i=1,\dots,4$, and
\begin{align*}
  p_1 = \lambda_1 + C x_3,\quad p_2 = \lambda_2 + C x_4,\quad p_{3} = \lambda_3, \quad p_{4} = \lambda_4.
\end{align*}
Then in $(q,p)$ coordinates, the Hamiltonian function takes the form:
\begin{align}\label{eq:Hc}
  H_C(q, p)=\frac{1}{2}\left(\frac{p_1^2}{\sigma}+p_2^2+\frac{p_3^2}{\sigma}+p_4^2-\frac{\tau(p_3 q_1-p_1 q_3+\sigma(p_4 q_2-p_2 q_4))^2}{\sigma(\sigma+(q_1^2+\sigma q_2^2+q_3^3+\sigma q_4^2)\tau)}\right).
\end{align}
The Poisson structure on the restricted manifold is given by:
\begin{align*}
  \{q_k,p_j\}&=\delta_{kj},\quad \{p_k,p_{j+n}\}=\delta_{kj}C,\quad
   \{q_k, q_j\}=\{p_k, p_j\}=0.
\end{align*}
 Hence, this system describes the motion of a charged particle in the constant magnetic field.

The Hamilton's system is:
\begin{gather}\label{eq:HameqR}
\begin{aligned}
    \dot{q_1}&=\{q_1, H_C\} = \frac{1}{\sigma}(p_1+W q_3), & \dot{p_1}&=\{p_1, H_C\} =(C+W)\dot{q_3},\\
    \dot{q_2}&=\{q_2, H_C\} =(p_2+W q_4),& \dot{p_2}&=\{p_2, H_C\} =(C+W)\dot{q_4},\\
    \dot{q_3}&=\{q_3, H_C\} =  \frac{1}{\sigma}(p_3-W q_1),& \dot{p_3}&=\{p_3, H_C\} =-(C+W)\dot{q_1},\\
    \dot{q_4}&=\{q_4, H_C\} =(p_4-W q_2), & \dot{p_4}&=\{p_4, H_C\} =-(C+W)\dot{q_2},
    \end{aligned}
    \end{gather}
where 
$ W=\dfrac{\tau(p_3 q_1-p_1 q_3+\sigma(p_4 q_2-p_2 q_4)}{\sigma+(q_1^2+\sigma q_2^2+q_3^3+\sigma q_4^2)\tau}$.

It is not hard to check that for $i=1,2$:
\begin{align*}
    \{q_{i+2}p_i-q_i p_{i+2},H_C\}&={\sigma^{i-2}}(p_{i+2}-W q_i)p_i+(C+W)q_{i+2}\dot{q}_{i+2}-{\sigma^{i-2}}(p_i+W q_{i+2})p_{i+2}+(C+W)q_i\dot{q}_i\\
    &={\sigma^{i-2}}(C+W)(q_i(p_i+W q_{i+2})+q_{i+2}(p_{i+2}-W q_i))-{\sigma^{i-2}}W(q_i p_i+q_{i+2} p_{i+2})\\
    &=C{\sigma^{i-2}}(q_i p_i+q_{i+2}p_{i+2})=\frac{C}{2}\{q_i^2+q_{i+2}^2,H_C\}.
\end{align*}
Hence,
\begin{align}\label{eq:lin}
    I_1 = C(q_1^2+q_3^2)-2(q_3 p_1-q_1 p_3),\quad
    I_2 = C(q_2^2+q_4^2)-2(q_4 p_2-q_2 p_4)
\end{align}
are the first integrals that are linear in the momenta. 

Also, we have two first integrals quadratic in the momenta. We have that the following holds:
\begin{align*}
\{\sum_{i=1}^2\sigma^{2-i}(p_i^2+p_{i+2}^2), H_C\}&=2\sum_{i=1}^2\sigma^{2-i}(C+W)(\dot{q}_{i+2}p_i-\dot{q}_i p_{i+2})\\
& =2(C+W)\sum_{i=1}^2 ((p_{i+2}-Wq_i)p_i-(p_i+W q_{i+2}) p_{i+2})\\
&=-2 W(C+W)\sum_{i=1}^2(q_i p_i+q_{i+2} p_{i+2})=-2 W(C+W)\sum_{k=1}^4 q_k p_k,\\
\{\sum_{k=1}^4 q_kp_k,H_C\}&=\sum_{i=1}^2\dot{q}_i(p_i-(C+W)q_{i+2})+\dot{q}_{i+2}(p_{i+2}+(C+W)q_i)\\
&=\sum_{i=1}^2\sigma^{i-2}((p_i+W q_{i+2})(p_i{-}(C+W)q_{i+2})+(p_{i+2}-W q_i)(p_{i+2}{+}(C+W)q_i))\\
&=\sum_{i=1}^2\sigma^{i-2}((p_i^2+p_{i+2}^2){+}C(p_{i+2} q_i-p_i q_{i+2}){-}W(C+W)(q_i^2+q_{i+2}^2))
\end{align*}
Now:
\begin{align*}
&\{\sigma(p_1^2+p_3^2)+p_2^2+p_4^2 + \tau(q_1p_1+q_2p_2+q_3p_3+q_4p_4)^2,H_C\}\\
=&2\left(\tau\sum_{i=1}^2\sigma^{i-2}((p_i^2+p_{i+2}^2){+}C(p_{i+2} q_i-p_i q_{i+2}){-}W(C+W)(1+\tau\sum_{j=1}^2\sigma^{j-2}(q_j^2+q_{j+2}^2))\right)\sum_{k=1}^4 q_k p_k\\
=&2\tau\left(\sum_{i=1}^2\sigma^{i-2}(p_i^2+p_{i+2}^2){-}W(p_{i+2} q_i-p_i q_{i+2}))\right)\sum_{k=1}^4 q_k p_k\\
=&4\tau H_C\sum_{k=1}^4 q_k p_k=2\tau\{\sum_{i=1}^2\sigma^{2-i}(q_i^2+q_{i+2}^2)H_C,H_C\}
\end{align*}
Therefore, the first quadratic integral is:
\begin{align}\label{eq:quad}
   I_3=\sigma(p_1^2+p_3^2)+p_2^2+p_4^2 + \tau(q_1p_1+q_2p_2+q_3p_3+q_4p_4)^2 - 2\tau(\sigma(q_1^2+q_3^2)+q_2^2+q_4^2) H_C,
\end{align}
 while the second is, naturally, the Hamiltonian $I_4=H_C$. 
 These integrals are functionally independent and commuting, hence the system is completely integrable.
\end{proof}

\begin{corollary}\label{cor:HLR}
    Hamiltonian system~\eqref{eq:HamReq} is integrable in Liouville sense via first integrals:
    \begin{align}\label{eq:intLR}
    \begin{aligned}
       I_0&=H,\quad I_1=\lambda_5,\\
       I_2 &= C(x_1^2-x_3^2)-2(x_3 \lambda_1-x_1 \lambda_3),\quad
    I_3 = C(x_2^2-x_4^2)-2(x_4 \lambda_2-x_2 \lambda_4),\\ 
    I_4&=\sigma((\lambda_1+C x_3)^2+\lambda_3^2)+(\lambda_2+Cx_4)^2+\lambda_4^2 \\
    &\phantom{=}+ \tau(x_1(\lambda_1+C x_3)+x_2(\lambda_2+Cx_4)+x_3\lambda_3+x_4\lambda_4)^2- 2\tau(\sigma(x_1^2+x_3^2)+x_2^2+x_4^2) H,
    \end{aligned}
    \end{align}
     with $H$ being the Hamiltonian function~\eqref{eq:HamR}.
\end{corollary}

\begin{remark}
    Note that LL system is superintegrable since its first integrals $I_0$, $I_1$, $I_2$, $I_5$ and $I_6$ from~\eqref{eq:intLL} are pairwise in involution. Moreover, there are two additional integrals $I_3$ and $I_4$ satisfying $\{I_1, I_3\}=\{I_2, I_4\}=I_0$ and $H=\frac{1}{2\sigma}(I_5+\sigma I_6)$. In contrast, LR system is only Liouville integrable since it has exactly five functionally independent first integrals~\eqref{eq:intLR} in involution. 
\end{remark}

A classification of the geodesic curves of the sub-Riemannian LR system on $\HH_5$ for the standard metric ($\sigma=\tau=1$) was presented in detail in~\cite{PAV}. The trajectories can be classified into two categories: those that lie outside the hypersphere $|q|^2=r_0$ (Figure~\ref{fig:geodH5} \textit{left}), and those that remain within the region bounded by the two hyperspheres $|q|^2= r_1$ and $|q|^2= r_2$ (Figure~\ref{fig:geodH5} \textit{right}). However, since their explicit parametrization involves elliptic integrals, we omit it here. To visualize objects in five-dimensional real space, we interpret the coordinate pairs $(x_1, x_3)$ and $(x_2, x_4)$ as complex numbers $z_1$ and $z_2$, respectively, and represent their moduli along the first two coordinates of the plot, while retaining the fifth coordinate as the third axis. The resulting three-dimensional plot therefore has coordinates $(\|z_1\|,\|z_2\|, z)$.

\begin{figure}[H]
    \centering
    \includegraphics[height=6cm]{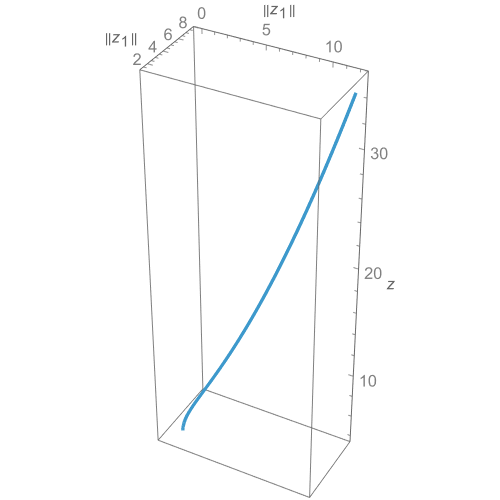}\qquad
    \includegraphics[height=6cm]{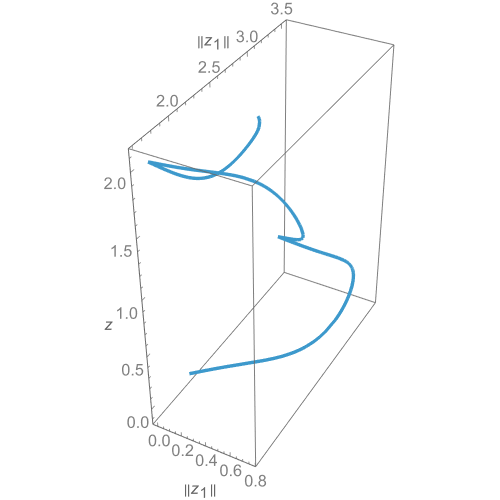}
    \caption{Examples of sub-Riemannian geodesic curves of LR system on $\HH_5$ for the standard metric.}
    \label{fig:geodH5}
\end{figure}

Finally, we present figure illustrating the difference between LL or LR curves for the restricted system~\eqref{eq:HameqR}.

\begin{figure}[H]
    \centering
    \includegraphics[width=0.5\linewidth]{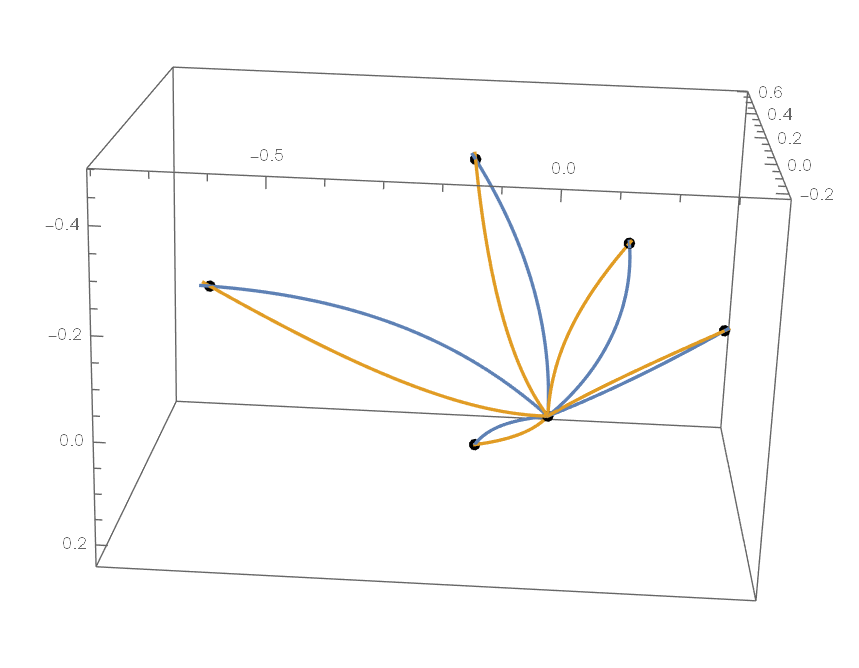}
    \caption{Stereographic projection of the trajectories for the restricted LL system (\textit{blue}) and LR system (\textit{orange}) connecting the origin to the same five randomly selected points.}
    \label{fig:geodH5}
\end{figure}

\begin{thank}
This research  is partially supported by the Ministry of  Science, Technological Development and Innovation, Republic of Serbia, through the project 451-03-33/2026-03/200104.
\end{thank}

\bibliographystyle{amsplain}

\end{document}